\documentclass[11pt]{article}
\usepackage{amsthm, amsmath, amssymb, amsfonts, url, booktabs, tikz, setspace, fancyhdr}
\usepackage[hidelinks]{hyperref}
\usepackage[margin = 1in]{geometry}
\usepackage{hyperref, enumerate}

\newtheorem{theorem}{Theorem}[section]
\newtheorem{proposition}[theorem]{Proposition}
\newtheorem{lemma}[theorem]{Lemma}
\newtheorem{corollary}[theorem]{Corollary}
\newtheorem{conjecture}[theorem]{Conjecture}

\theoremstyle{definition}

\newtheorem{example}[theorem]{Example}

\theoremstyle{remark}

\newcommand{\norm}[1]{\left\lVert#1\right\rVert}
\newcommand{\ang}[1]{\left\langle #1 \right\rangle}

\newcommand{\overbar}[1]{\mkern 1.5mu\overline{\mkern-1.5mu#1\mkern-1.5mu}\mkern 1.5mu}

\newcommand{\RR}{\mathbb{R}}

\newcommand{\HH}{\mathrm{H}}
\newcommand{\cP}{\mathcal{P}}
\newcommand{\FF}{\mathcal{F}}
\newcommand{\TT}{\mathcal{T}}

\newcommand{\A}{\mathcal{A}}

\newcommand{\C}{\mathcal{C}}

\newcommand{\J}{\mathcal{J}}

\newcommand{\ee}{\mathbf{e}}
\newcommand{\db}{\mathrm{db}}
\newcommand{\Hom}{\mathrm{Hom}}
\newcommand{\M}{\mathcal{M}}
\newcommand{\ww}{\mathbf{w}}

\usetikzlibrary{trees}
\tikzstyle{p}+=[fill=black, circle, minimum width = 1pt, inner sep =
1pt]
\tikzstyle{w}+=[fill=white, draw, circle, minimum width = 1pt, inner sep =
1.5pt]

\newcommand{\red}[1]{{\color{red}#1}}

\begin{document}

\title{Finite reflection groups and graph norms}

\author{David Conlon\thanks{Mathematical Institute, Oxford OX2 6GG, United
Kingdom.
E-mail: {\tt
david.conlon@maths.ox.ac.uk}. Research supported by a Royal Society University Research Fellowship and by ERC Starting Grant 676632.} \and Joonkyung Lee\thanks{Mathematical Institute, Oxford OX2 6GG, United
Kingdom.
E-mail: {\tt
joonkyung.lee@maths.ox.ac.uk}. Supported by the ILJU Foundation of Education and Culture.}}

\date{}

\maketitle

\begin{abstract}
Given a graph $H$ on vertex set $\{1,2,\cdots, n\}$ and a function $f:[0,1]^2\rightarrow \RR$, define 
\begin{align*}
\norm{f}_{H}:=\left\vert\int \prod_{ij\in E(H)}f(x_i,x_j)d\mu^{|V(H)|}\right\vert^{1/|E(H)|},
\end{align*}
where $\mu$ is the Lebesgue measure on $[0,1]$. We say that $H$ is norming if $\norm{\cdot}_H$ is a semi-norm.
A similar notion $\norm{\cdot}_{r(H)}$ is defined by $\norm{f}_{r(H)}:=\norm{|f|}_{H}$
and $H$ is said to be weakly norming if $\norm{\cdot}_{r(H)}$ is a norm.
Classical results show that weakly norming graphs are necessarily bipartite. In the other direction, Hatami showed that even cycles, complete bipartite graphs, and hypercubes are all weakly norming. 
We demonstrate that any graph whose edges percolate in an appropriate way under the action of a certain natural family of automorphisms is weakly norming. 
This result includes all previously known examples of weakly norming graphs, but also allows us to identify a much broader class arising from finite reflection groups. We include several applications of our results. In particular, we define and compare a number of generalisations of Gowers' octahedral norms and we prove some new instances of Sidorenko's conjecture.
\end{abstract}

\section{Introduction}
Let $H$ be a graph on vertex set $\{1,2,\cdots,n\}$ and 
$f:[0,1]^2\rightarrow \RR$ be a bounded Lebesgue measurable function. 
Consider the integral
	\begin{align}\label{eqn:t_H}
		\int \prod_{ij\in E(H)}f(x_i,x_j)d\mu^{|V(H)|},
	\end{align}
where $\mu$ is the Lebesgue measure on $[0,1]$.
If we choose $f$ so as to model the adjacency matrix of a graph $G$, the integral above corresponds to the homomorphism density $t_H(G)$ of $H$ in $G$, which plays a central role in extremal graph theory. Similar expressions also appear naturally in other areas, particularly in statistical physics. 

Our concern in this paper will be with the natural question, proposed by Lov\'asz \cite{L12}, of determining when the integral \eqref{eqn:t_H} defines a (semi-)norm. Formally, we say that a graph $H$ is \emph{norming} if the functional defined by
 	\begin{align}\label{eqn:H-norm}
		\norm{f}_{H}:=\left\vert\int \prod_{ij\in E(H)}f(x_i,x_j)d\mu^{|V(H)|}\right\vert^{1/|E(H)|}
	\end{align}
is a semi-norm, and \emph{weakly norming} if
	\begin{align}\label{eqn:abs_H-norm}
		\norm{f}_{r(H)}:=\left(\int \prod_{ij\in E(H)}|f(x_i,x_j)|d\mu^{|V(H)|}\right)^{1/|E(H)|}
	\end{align}
is a norm. As one might expect from the name, 
it is easy to check that every norming graph is weakly norming.
While our focus will usually be on weakly norming graphs,
we will often discuss how analogous results can be derived for norming graphs under an extra technical condition.

The first in-depth study of (weakly) norming graphs was undertaken by Hatami \cite{H10}. A moment's thought shows that $H$ is necessarily bipartite whenever it is weakly norming, because otherwise it can be zero for strictly positive functions.
In \cite{H10}, Hatami showed that the $n$-dimensional hypercube $Q_n$ and the complete bipartite graph $K_{m,n}$ are weakly norming.
He also observed that the functionals $\norm{\cdot}_{C_{2k}}$ correspond to the classical Schatten--von Neumann norms
and, hence, even cycles are norming. In addition, Lov\'{a}sz~\cite{L12} showed that the complete bipartite graph $K_{n,n}$ minus a perfect matching is weakly norming.

We generalise these results, finding a much larger class of (weakly) norming graphs that includes all of the known examples.
To give an indication of our results, suppose that $k$ and $r$ are integers with $k\leq r$ and $\mathcal{P}$ is a polytope. Consider the bipartite graph between $k$-faces and $r$-faces of $\mathcal{P}$ indicating their incidence. That is, we place an edge between a $k$-face and an $r$-face if one contains the other.
We call this graph the $(k,r)$-incidence graph of the polytope $\mathcal{P}$. We then have the following theorem.

\begin{theorem}\label{thm:polytope_norming}
$(k,r)$-incidence graphs of regular polytopes are weakly norming.
\end{theorem}

For example, in an $n$-dimensional simplex, the $k$-faces and $r$-faces naturally correspond to $(k+1)$-element and $(r+1)$-element subsets of $[n]$. 
Therefore, the $(k,r)$-incidence graph of an $n$-simplex is the inclusion graph between $(k+1)$-sets and $(r+1)$-sets.
In particular, the $(0,1)$-incidence graph is the 1-subdivision of $K_n$,
the $(0,n-2)$-incidence graph is $K_{n,n}$ minus a perfect matching, and
the $(0,n-1)$-incidence graph is the star $K_{1,n}$, which by tensor powering shows that $K_{m,n}$ is also weakly norming.
Even cycles $C_{2k}$ are the $(0,1)$-incidence graphs of regular $k$-gons and thus are weakly norming.
More generally, by considering the $(0,1)$-incidence graph of any regular polytope, such as hypercubes or the icosahedron,
we see that their 1-subdivisions are weakly norming. 

We prove Theorem \ref{thm:polytope_norming} as a corollary of a more general result showing how weakly norming graphs arise from finite reflection groups.
A \emph{finite reflection group} is
a finite subgroup of $\mathbf{GL}(n,\RR)$ generated by a set of reflections across hyperplanes passing through the origin.
Those readers who are not familiar with these groups
may temporarily assume that, given a finite reflection group $W$, 
there exists a distinguished set of generators $S$ known as simple reflections (for further details,
see Section~\ref{sec:Prelim}).
Fixing a generating set of simple reflections $S$ in a finite reflection group $W$, let $S_1$ and $S_2$ be subsets of $S$ and 
let $W_1$ and $W_2$ be the subgroups of $W$ generated by $S_1$ and $S_2$, respectively.
Consider the bipartite graph between the (left-)cosets of $W_1$ and $W_2$,
where $wW_1$ and $wW_2$ are adjacent for every $w\in W$.
We call this graph the \emph{$(S_1,S_2;S,W)$-reflection graph}
and say that a graph $H$ is a \emph{reflection graph} if it is isomorphic to an $(S_1,S_2;S,W)$-reflection graph for some suitable choice of parameters.
With these definitions, we may now state our main result.

\begin{theorem}\label{thm:parabolic}
Reflection graphs are weakly norming.
\end{theorem}

This class includes the $(k, r)$-incidence graphs of regular polytopes, but also provides other simple examples, 
such as the hypercube (rather than its subdivision) and the graph obtained by replacing each edge of an octahedron with a cycle of length $4$ (see Example~\ref{ex:norming_B3} for a more formal description). 
It also opens the door to more exotic examples, coming from the exceptional reflection groups $E_6$, $E_7$, and $E_8$. An analogous result also holds for norming graphs, though under a slightly stronger condition.

\begin{theorem}\label{thm:norming}
Let $W$ be a finite reflection group and let $S$ be a generating set of simple reflections.
Then, for any disjoint subsets $S_1$ and $S_2$ of $S$, the $(S_1,S_2;S,W)$-reflection graph is norming.
\end{theorem}

When proving that $\norm{\cdot}_{r(H)}$ is a norm, all of the difficulties lie in proving the triangle inequality. 
Hatami's work in \cite{H10} started from the observation that a H\"older-like inequality is equivalent to the triangle inequality for $\norm{\cdot}_{r(H)}$.
To state his condition, we have to introduce some notation generalising \eqref{eqn:H-norm} and \eqref{eqn:abs_H-norm}.
Let $m=|E(H)|$ and let $\chi:E(G)\rightarrow [m]$ be a (not necessarily proper) edge colouring  of $H$.
Consider a family $\FF=\{f_1,f_2,\cdots,f_m\}$ of bounded measurable functions on $[0,1]^2$
and define 
	\begin{align*}
		\ang{\FF;\chi}_{H}:=\int \prod_{e=ij\in E(H)}f_{\chi(e)}(x_i,x_j)d\mu^{|V(H)|}.
	\end{align*}
Note that if $f_i=|f|$ for all $i=1,2,\cdots,m$, then $\ang{\FF;\chi}_H=\norm{f}_{r(H)}^{|E(H)|}$, while if $f_i = f$, then
$|\ang{\FF;\chi}_{H}|=\norm{f}_H^{|E(H)|}$.
Hatami's result now says that the triangle inequality holds for $\norm{\cdot}_{r(H)}$ if and only if
	\begin{align}\label{eqn:rainbow}
		\ang{\FF;\chi}_H\leq \prod_{e \in E(H)}\norm{f_{\chi(e)}}_{r(H)}
	\end{align}
for all choices of $\FF$ and $\chi$. 
Furthermore, $\norm{\cdot}_H$ is a semi-norm if and only if the analogous inequality obtained by replacing $\norm{f_{\chi(e)}}_{r(H)}$ with $\norm{f_{\chi(e)}}_H$ holds.

We will think of \eqref{eqn:rainbow} in the following terms:
regard the functions $f_1,\cdots,f_m$ as $m$ distinct graphs on the same vertex set\footnote
{When the $f_i$ are non-negative, as they will be when studying weakly norming graphs, there are large graphs approximating $f_{1}, \cdots
,f_{m}$ by the limit theory of dense graphs \cite{LSz06}.} 
and imagine each edge of $f_i$ has the colour $i$.
Then $\ang{\FF;\chi}_H$ is the number of (homomorphic) copies of $H$ which are coloured according to $\chi$,
i.e., each edge $e\in E(H)$ receives the colour $\chi(e)$.
In particular, if $\chi$ is a one-to-one map then $\ang{\FF;\chi}_H$ counts the number of `rainbow' copies of $H$, 
while $\norm{f_i}_{r(H)}^{|E(H)|}$ counts the number of monochromatic copies of $H$ in colour $i$.
Thus, \eqref{eqn:rainbow} is equivalent to the statement that the number of rainbow copies of $H$ is bounded 
above by the geometric mean of the number of monochromatic copies in each colour.

The proof that~\eqref{eqn:rainbow} holds for reflection graphs has two steps. In the first step, discussed in Section~\ref{sec:CS}, we show that any graph whose edges percolate in an appropriate way under the action of a certain natural family of automorphisms is weakly norming. This statement, Theorem~\ref{cor:colouring_graphs}, is more general than Theorem~\ref{thm:parabolic}, and may be of independent interest, but the resulting condition needs to be verified by hand for any particular graph or class of graphs. Accordingly, the second step in our proof, discussed in Section~\ref{sec:euc}, is to find a general argument that verifies this condition for all reflection graphs. It is here that we use results from the theory of finite reflection groups.

Suppose now that $H$ is a weakly norming graph, $f$ is a bounded measurable function on $[0,1]^2$, and $e^*$ is an edge of $H$. If we put $f_1= |f|$, $f_2=f_3=\cdots=f_m=1$, and make $\chi$ one-to-one with $\chi(e^*)=1$, 
then $\ang{\FF;\chi}_H=\norm{f}_{r(K_2)}$, where $K_2$ is just a single edge, so~\eqref{eqn:rainbow} implies that
\begin{align}\label{eqn:Sidorenko}
	\norm{f}_{r(K_2)}\leq \norm{f}_{r(H)}.
\end{align}
That is, when $H$ is weakly norming, $H$ satisfies Sidorenko's conjecture~\cite{Sid93, Sid92}, which says exactly that for any bipartite graph $H$ and any bounded measurable function $f$, an inequality of the form~\eqref{eqn:Sidorenko} holds.
Sidorenko's conjecture is one of the major open problems in extremal graph theory, 
and there has been much recent work \cite{CFS10, CKLL15, KLL14, LSz12, Sz15} verifying the conjecture for a widening class of graphs.
As noted above, all of the weakly norming graphs found in this paper also satisfy Sidorenko's conjecture.
However, this is not the only application of our results to Sidorenko's conjecture. By applying the entropy techniques developed in~\cite{LSz12, KLL14, Sz15, CKLL15}, we will prove that weakly norming graphs 
can also be used as building blocks for constructing new graphs that satisfy the conjecture. We refer the reader to Section~\ref{sec:Sid} for more details. 

Finally, also in Section~\ref{sec:genapp}, we will discuss generalisations of our results to hypergraphs, with the main result being that a suitably defined family of reflection hypergraphs are weakly norming. We will then show that every norm defined in this manner is equivalent, in some well-defined sense, to a corresponding cut-norm and then to an appropriate variant of the octahedral norms introduced by Gowers~\cite{G06, G07} in his work on hypergraph regularity. Our methods also allow us to compare the relative strengths of these variants. These results generalise earlier work by Gowers~\cite{G06, G07} and by Conlon, H\`an, Person, and Schacht~\cite{CHPS12} on equivalences between cut-norms and octahedral norms and the notions of quasirandomness they define.

\section{A motivating example}\label{sec:C6}

It is already non-trivial to show that \eqref{eqn:Sidorenko} holds, even for graphs as simple as paths~\cite{BR65} or trees~\cite{Sid93, Sid92}, so it should 
not be surprising that it is more difficult to prove the strictly stronger inequality \eqref{eqn:rainbow}.
To motivate what follows, we will prove that $C_6$, the cycle of length six, is norming without invoking spectral graph theory or Schatten--von Neumann norms.
Let $H=C_6$ be the graph with vertex set $\{1, 2, \cdots,6\}$, where $i$ and $i+1$ are adjacent for all $i$ (with addition taken modulo $6$),
and let $\FF=\{f_1,f_2,\cdots,f_6\}$ be a family of six functions, each taking two variables.

To show that \eqref{eqn:rainbow} holds, we may assume that $\chi$ is a one-to-one map, that is, a rainbow colouring, as it is in the most general case.
Without loss of generality, put $\chi(e_i)=i$ if $e_i$ is the edge between $i$ and $i+1$ modulo 6.
Define functions $g$ and $h$  by
	\begin{align*}
		g(x_1,x_4)&=\int f_{1}(x_1,x_2)f_{2}(x_2,x_3)f_{3}(x_3,x_4)~dx_2 dx_3 \text{, and}\\
		h(x_1,x_4)&=\int f_{4}(x_4,x_5)f_{5}(x_5,x_6)f_{6}(x_6,x_1)~dx_5 dx_6,
	\end{align*} 
so that $\int gh~dx_1 dx_4= \ang{\FF;\chi}_{C_6}$.
The function $g$ gives the (normalised) count of rainbow walks from $x_1$ to $x_4$ which are coloured, in order, with the colours $1$, $2$, and $3$. The function $h$ can be similarly interpreted, but using the colours $6$, $5$, and $4$. A simple application of the Cauchy--Schwarz inequality gives
	\begin{align}\label{eqn:basic_CS}
		\ang{\FF;\chi}^2_{C_6}\leq \int g^2 \int h^2,
	\end{align}
where here and throughout the paper we suppress the variables of integration if they are clear from context.
Observe now that
	\begin{align*}
		\int g^2 = \ang{\FF;\chi_L}_{C_6},	
	\end{align*}
where $\chi_L$ is the colouring obtained by doubling the `left half' of the rainbow cycle, represented by the vector $(1,2,3,3,2,1)$.
A similar equation also holds for $h$, namely, $\int h^2 = \ang{\FF;\chi_R}_{C_6}$,
where $\chi_R=(6,5,4,4,5,6)$ is the colouring obtained by doubling the `right half' of the rainbow cycle. We have therefore bounded  
$\ang{\FF;\chi}_{C_6}$ from above by the geometric mean of two functions of the same form, but simpler in the sense that they both contain fewer colours.
Repeating this procedure twice more, by first applying the Cauchy--Schwarz inequality to each of $\ang{\FF;\chi_L}_{C_6}$ and $\ang{\FF;\chi_R}_{C_6}$ with respect to the variables $(x_2, x_5)$ and then applying the Cauchy--Schwarz inequality to each of the four resulting terms with respect to the variables $(x_3, x_6)$, we find that 
\begin{align*}
	\ang{\FF;\chi}_{C_6}\leq\prod_{i=1}^8 \ang{\FF;\chi_i}^{1/8}_{C_6},
\end{align*}
where 
\begin{align}\label{eqn:3-step_C6}
&\chi_1=(1,1,1,1,1,1), \ \chi_2=(2,1,1,2,2,2), \ \chi_3=(2,2,2,2,3,3), \ \chi_4=(3,3,3,3,3,3),\nonumber\\
&\chi_5=(6,6,6,6,6,6), \ \chi_6=(5,6,6,5,5,5), \ \chi_7=(5,5,5,5,4,4), \ \chi_8=(4,4,4,4,4,4).
\end{align}
Recall that $\ang{\FF;\chi_1}_{C_6}= \norm{f_1}_{C_6}^6$ 
and hence it gives one of the terms in the desired upper bound,
though with an incorrect exponent. 

We now iterate this $3$-step process of applying the Cauchy--Schwarz inequality along different vertex cuts, noting that the more
often we repeat the process the more 
monochromatic colourings appear in our upper bound.
In particular, as in~\eqref{eqn:3-step_C6}, after each repetition, the sum of the exponents of the
non-monochromatic forms reduces by at least a half.
Hence, if we iterate $k$ times, we get the bound
\begin{align*}
	\ang{\FF;\chi}_{C_6}\leq \prod_{i=1}^{2^{3k}}\ang{\FF;\chi_{i,k}}_{C_6}^{1/2^{3k}},
\end{align*}
where at least a $1 - 2^{-k}$ proportion of the $2^{3k}$ colourings $\chi_{i,k}$ are monochromatic. We may rewrite this inequality as
\begin{align*}
	\ang{\FF;\chi}_{C_6}\leq \prod_{i=1}^{6^6}\ang{\FF;\chi'_{i}}_{C_6}^{\alpha_{i,k}},
\end{align*}
where each $\chi'_{i}$ represents one of the $6^6$ possible edge-colourings of $C_6$ with $6$ colours. 
In particular, we assume that
$\chi'_{1}, \chi'_2, \cdots, \chi'_{6}$ represent the monochromatic colourings in $1, 2, \cdots, 6$, respectively. By the argument above, $\alpha_{1,k} + \dots + \alpha_{6,k} \geq 1 - 2^{-k}$ and it is also easy to see that $\alpha_{i,k}$ is non-decreasing in $k$ for all $i = 1, 2, \cdots, 6$. Therefore, taking the limit as $k$ tends to infinity, we have
\begin{align*}
	\ang{\FF;\chi}_{C_6}\leq \prod_{i=1}^{6}\ang{\FF;\chi'_{i}}_{C_6}^{\alpha_{i}},
\end{align*}
where $\alpha_1 + \alpha_2 + \cdots + \alpha_6 = 1$. If $\alpha_1, \alpha_2, \cdots, \alpha_6$ were equal, this would be the desired inequality. If they are not equal, we note that an analogous procedure, but applying the Cauchy--Schwarz inequalities first with $(x_2, x_5)$, then with $(x_3, x_6)$, and finally with $(x_1, x_4)$, allows one to prove the inequality
\begin{align*}
	\ang{\FF;\chi}_{C_6}\leq \prod_{i=1}^{6}\ang{\FF;\chi'_{i}}_{C_6}^{\beta_{i}},
\end{align*}
where $\beta_i = \alpha_{i-1}$ (and addition is again taken modulo $6$). Repeating the same idea four more times, we can cyclically permute the exponents in the inequality to all six possible positions. Taking the product of these six inequalities then completes the proof.

\section{Cut involutions and Cauchy--Schwarz inequalities} \label{sec:CS}

There were three steps to the proof given in the previous section: 
firstly, we showed how to apply the Cauchy--Schwarz inequality along certain vertex cuts; then we showed that
a monochromatic edge-colouring of $H$ can be obtained
through a sequence of such Cauchy--Schwarz inequalities;
finally, we used a limiting argument and the edge transitivity of $C_6$ to complete the proof. In this section,
we generalise these arguments. 
The first step will be generalised through the use of a natural class of graph automorphisms which we call cut involutions, 
while an appropriate generalisation of the last step is relatively straightforward. 
Generalising the second step, that is, finding a monochromatic edge-colouring, proves more difficult,
so in this section we reduce it to a simpler question which we will resolve in Section~\ref{sec:euc} for the special case of reflection graphs.

To generalise the Cauchy--Schwarz inequality \eqref{eqn:basic_CS},
we need to find a vertex cut which splits the graph into two isomorphic pieces,
which we refer to, informally, as the left and right halves.
To be more precise, given a connected graph $H$, say that an automorphism $\phi$ of $H$ is a \emph{cut involution} if it is an involution, i.e., $\phi=\phi^{-1}$, 
and the fixed point set $F_\phi = \{v\in V(H):\phi(v)=v\}$ is a vertex cut of $H$.
Let $W_H$ be the subgroup of the automorphism group of $H$ generated by the cut involutions and call it the \emph{cut involution group} of $H$.

Since it is possible for $H\setminus F_{\phi}$ to have more than two components,
the left and right halves of a cut involution may not be well-defined.
In reflection graphs, we can avoid this ambiguity by focusing on certain special cut involutions.
For now, we will simply assume that the left and right halves
are somehow fixed for each cut involution $\phi$.
That is, we will assume that there is a fixed pair of disjoint subsets $L_{\phi}$ and $R_{\phi}$ of $V(H)\setminus F_\phi$
that are unions of connected components in $H\setminus F_\phi$ and which are mapped to each other under $\phi$.  

Define the \emph{left-folding map} $\phi^+:V(H)\rightarrow V(H)$ of a cut involution $\phi$ by
\begin{align*}
	\phi^+(v)=
	\begin{cases}
		\phi(v) ~&\text{ if }v\in R_\phi\\
		v &\text{ if }v\in L_\phi\cup F_\phi,
	\end{cases}
\end{align*}
and, similarly, define the \emph{right-folding map} $\phi^-$ of $\phi$ by swapping the roles of $L_\phi$ and $R_\phi$.
Let $m=|E(H)|$, $\FF=\{f_1, f_2, \cdots, f_m\}$ be a family of non-negative bounded measurable functions on $[0,1]^2$, and $\chi$ be an edge-colouring of $H$ with colours $\{1,2,\cdots,m\}$.
Then the Cauchy--Schwarz inequality \eqref{eqn:basic_CS} can be generalised in terms of a cut involution $\phi$ by
\begin{align}\label{eqn:CS}
	\ang{\FF;\chi}_H\leq 
	\ang{\FF;\chi\circ\phi^+}_{H}^{1/2}\ang{\FF;\chi\circ\phi^-}_{H}^{1/2}.
\end{align}
Here $\chi\circ\phi^+$ is the colouring of $H$ where the colour of the edges in $R_\phi\cup F_\phi$ copies the corresponding edges in $L_\phi\cup F_\phi$. Similarly, $\chi\circ\phi^-$ is the colouring where the colour of the edges in $L_\phi\cup F_\phi$ copies the corresponding edges in $R_\phi\cup F_\phi$.

Note that in \eqref{eqn:CS}, the non-negativity of the functions in $\FF$ is used only if there is an edge inside $F_{\phi}$.
More precisely, if the fixed point set $F_\phi$ contains edges, then we are using a Cauchy--Schwarz inequality of the form
\begin{align*}
	\int hg_1g_2 =\int (h^{1/2}g_1)(h^{1/2}g_2) \leq \int hg_1^2 \int hg_2^2, 
\end{align*}
where $h=\prod_{(u,v)\in H[F_\phi]}f_{\chi(uv)}(x_u,x_v)$.
Conversely, if there is no edge of $H$ fixed by $\phi$, then the non-negativity of functions in $\FF$ is unnecessary.
This observation leads naturally to the definition of a \emph{stable involution}. This is a cut involution $\phi$ such that the fixed point set $F_{\phi}$ contains no edge of $H$, that is, it is an independent (or stable) set in $H$ as well as a vertex cut. The subgroup of the automorphism group of $H$ generated by the stable involutions is then called the \emph{stable involution group} of $H$ and is denoted $S_H$.
In what follows, we will show that graphs whose edges percolate under the action of the cut involution group are weakly norming. Using the simple observation above, it is possible to prove analogous results for the stronger norming property by substituting the stable involution group for the cut involution group. We will not always be explicit about this below, choosing to focus instead on cut involutions and the weakly norming property. Nevertheless, it is worth bearing in mind.

Consider now a sequence $\phi_1,\phi_2,\cdots, \phi_d$ of cut involutions
and suppose that we wish to apply a sequence of Cauchy--Schwarz inequalities of the form 
\eqref{eqn:CS}, first with $\phi=\phi_1$, then with $\phi = \phi_2$, and so on. 
Let $\TT$ be the rooted binary tree of depth $d$ encoding 
which colourings have been obtained through such iterations: 
the vertices are labelled by $m$-edge-colourings of $H$, 
the root is labelled by the initial colouring $\chi_0$, and each vertex at depth $i<d$ labelled with $\chi$, say, 
has two children with labels $\chi\circ\phi_{i+1}^+$ and $\chi\circ\phi_{i+1}^-$.
We call this tree the \emph{Cauchy--Schwarz tree} associated with $(\phi_i)_{i=1}^{d}$ rooted at $\chi_0$. As here,
we will often abuse notation by identifying a vertex and its label.
\begin{example}\label{ex:CS-tree}
	The 3-step process described in the previous section
	corresponds to the following Cauchy--Schwarz tree of depth 3:
	\begin{center}
		\begin{tikzpicture}[level distance=1.5cm,
  level 1/.style={sibling distance=8.5cm},
  level 2/.style={sibling distance=4.25cm},
  level 3/.style={sibling distance=2.125cm} ]
  \node {(1,2,3,4,5,6)}
    child {node {(1,2,3,3,2,1)}
      child{node {(1,1,1,2,2,1)} 
      	child{node{(1,1,1,1,1,1)}}
      	child{node{(2,1,1,2,2,2)}}
      	}
      child {node {(2,2,3,3,3,3)}
      	child{node{(2,2,2,2,3,3)}}
      	child{node{(3,3,3,3,3,3)}}
      	}
    }
   child {node {(6,5,4,4,5,6)}
      child {node{(6,6,6,5,5,6)} 
      	child{node{(6,6,6,6,6,6)}}
      	child{node{(5,6,6,5,5,5)}}
      	}
      child{node {(5,5,4,4,4,4)}
      	child{node{(5,5,5,5,4,4)}}
      	child{node{(4,4,4,4,4,4)}}
      	}
    };
		\end{tikzpicture}
	\end{center}
\end{example}
Roughly speaking, a Cauchy--Schwarz tree shows how colours spread under applications of Cauchy--Schwarz inequalities of the form \eqref{eqn:CS}.
With this terminology, the following theorem generalises part of what was proven in the previous section.
\begin{theorem}\label{thm:percolate}
	Let $H$ be an edge-transitive graph with $m$ edges and $(\phi_i)_{i=1}^{d}$ be a finite sequence of cut involutions of $H$. 
	Let $\chi:E(G)\rightarrow [m]$ be a rainbow edge-colouring.
	If the Cauchy--Schwarz tree $\TT$ associated with $(\phi_i)_{i=1}^d$ rooted at $\chi$ contains 
	a leaf labelled with a monochromatic colouring $\chi'$,
	then $H$ is weakly norming.
\end{theorem}
\begin{proof}
	Consider the Cauchy--Schwarz tree $\TT$ associated with $(\phi_i)_{i=1}^{d}$ rooted at $\chi$,
	such that amongst its leaves $\chi_{1,d},\cdots,\chi_{2^d,d}$ there exists a leaf labelled with the monochromatic colouring $\chi'$.
	By applying \eqref{eqn:CS} repeatedly, we arrive at the upper bound
	\begin{align*}
		\ang{\FF;\chi}_H\leq 
		\prod_{i=1}^{2^d}\ang{\FF;\chi_{i,d}}_{H}^{1/2^d}.
	\end{align*}
	We now iterate the whole $d$-step process $k$ times.
	Consider the Cauchy--Schwarz tree $\TT'$ of depth $dk$ associated with the repeated sequence rooted at $\chi$.
	For $i=1,2,\cdots,m$, let $\chi_i'$ be the monochromatic colouring that only uses the colour $i$.
	Observe that all descendants of a vertex in $\TT'$ that is
	labelled with $\chi_i'$ are again labelled with $\chi_i'$.	
	Moreover, by our main assumption, every iteration of the $d$-step process makes 
	at least a $1/2^d$ proportion of the non-monochromatic leaves monochromatic.
	We may rearrange the upper bound obtained from $\TT'$ as
	\begin{align*}
		\ang{\FF;\chi}_H\leq 
		\prod_{i=1}^{m^m}\ang{\FF;\chi_i'}_{H}^{\alpha_{i,k}},
	\end{align*}
	where each $\chi'_{i}$ represents one of the $m^m$ possible edge-colourings of $H$ with $m$ colours.
	By the remarks above, $\alpha_{1,k}+\cdots+\alpha_{m,k}\geq 1-(1-1/2^d)^k$ and 
	each $\alpha_{i,k}$ with $1 \leq i \leq m$ is non-decreasing in $k$.
	Therefore, taking the limit as $k$ tends to infinity, we get
	\begin{align*}
		\ang{\FF;\chi}_H\leq 
		\prod_{i=1}^{m}\ang{\FF;\chi_i'}_{H}^{\alpha_{i}},
	\end{align*}
	where $\alpha_{1}+\cdots +\alpha_{m}=1$.
	It remains to prove that we may take $\alpha_{i}=1/m$.
	Observe that if $\psi$ is an automorphism of $H$ then
	$\ang{\FF;\chi}=\ang{\FF;\chi\circ\psi}$,
	since the colouring $\chi \circ \psi$ may be seen as the same colouring but with the vertices relabelled. 
	On the other hand, $\psi$ can be regarded as a permutation of the set of colours $\{1,2,\cdots,m\}$:
	if an edge $e\in E(H)$ receives colour $j=\chi(e)$ under the rainbow colouring $\chi$, 
	then it receives $i=\chi(\psi(e))$ under $\chi\circ\psi$.
	Repeating the same argument as above, 
	but with a Cauchy--Schwarz tree rooted at $\chi\circ\psi$, 
	we have the inequality
	\begin{align*}
		\ang{\FF;\chi}=\ang{\FF;\chi\circ\psi}\leq \prod_{i=1}^{m}\ang{\FF;\chi_i'}_{H}^{\beta_{i}},
	\end{align*}
	where $\beta_{i}=\alpha_{j}$ if $i=\chi \circ \psi \circ \chi^{-1} (j)$.
	Taking the geometric mean of these upper bounds
	over all automorphisms $\psi$, we get
	\begin{align*}
		\ang{\FF;\chi}\leq \prod_{i=1}^{m}\ang{\FF;\chi_i'}_{H}^{\gamma_{i}},
	\end{align*}
	where 
	\begin{align*}
		\gamma_{i}=\frac{1}{|\mathrm{Aut}(H)|}
		(\alpha_{1}|A_{1\rightarrow i}|+\alpha_{2}|A_{2\rightarrow i}|+\cdots+\alpha_{m}|A_{m\rightarrow i}|).
	\end{align*}
	Here $\mathrm{Aut}(H)$ is the group of all automorphisms of $H$
	and $A_{j\rightarrow i}$ is the set of automorphisms sending the edge with colour $j$ to the edge with colour $i$.
	Note that all $A_{j\rightarrow i}$ have the same size,
	as they are all cosets of the subgroup $A_{i\rightarrow i}$.
	Therefore, $\gamma_{i}=1/m$ for all $i$, as required.
\end{proof}

The theorem above shows that if we want to prove that a graph $H$ is weakly norming, it is
sufficient to show that there exists a Cauchy--Schwarz tree with a rainbow root and a monochromatic leaf.
Suppose now that $J$ is an edge subset of $H$ and $\phi$ a cut involution of $H$.
We define two edge sets $J^+(\phi)$ and $J^-(\phi)$ as follows:
\begin{align*}
	J^+(\phi)=\{e\in E(H):\phi^+(e)\in E(J)\}\text{ and }
	J^-(\phi)=\{e\in E(H):\phi^-(e)\in E(J)\}.
\end{align*}
That is, $J^+(\phi)$ is the graph formed by copying the edges of $J$ from the left half onto the right half. Similarly, $J^-(\phi)$ copies the edges from the right half onto the left half.
Let $J_0,J_1,J_2,\cdots$ be a sequence of edge subsets of $H$. We say that it is a \emph{folding sequence} in $H$ if, for each $i\geq 0$,
\begin{align*}
	J_{i+1}=J_i^+(\phi)\text{ or }J_{i+1}=J_i^-(\phi)
\end{align*}
for some cut involution $\phi$.
If a finite folding sequence $J_0,J_1,\cdots,J_N$ in a graph $H$ 
starts from a set $J_0$ consisting of a single edge and ends with $J_N=E(H)$,
then we call it a \emph{percolating sequence}.
With this terminology, we may rephrase Theorem \ref{thm:percolate}.

\begin{theorem}\label{cor:colouring_graphs}
	Suppose that $H$ is a graph which is edge-transitive under the cut involution group $W_H$.
	Then, if there exists a percolating sequence $J_0,J_1,\cdots,J_N$, $H$ is weakly norming.
\end{theorem}

We say that a percolating sequence is a \emph{strong percolating sequence} if every cut involution used in the sequence is a stable involution. The analogue of Theorem~\ref{cor:colouring_graphs} for the full norming property is then as follows.

\begin{theorem}\label{cor:strong_colouring_graphs}
	Suppose that $H$ is a graph which is edge-transitive under the stable involution group $S_H$.
	Then, if there exists a strong percolating sequence $J_0,J_1,\cdots,J_N$, $H$ is norming.
\end{theorem}

The key question now is whether it is possible to percolate over all edges starting from a single edge.
Although there is always enough flexibility when choosing between $J_i^+(\phi)$ and $J_i^-(\phi)$
to guarantee that the number of edges does not decrease, 
this is far from guaranteeing that a percolating sequence exists.
Indeed, as we update, we may lose as well as gain edges, making
it difficult to keep track of the changes.
However, for any known example of a weakly norming graph, 
one may easily check that there exists a percolating sequence.
The main purpose of the next section is to find a common generalisation for these ad hoc arguments.

\section{Euclidean embeddings of reflection graphs} \label{sec:euc}

\subsection{Preliminaries on finite reflection groups}\label{sec:Prelim}

In this subsection, we state some preliminary facts about reflection groups, 
focusing on the case of finite groups to make the discussion more concise.
For more details, we refer the reader to \cite{BB05} and \cite{H92}, while those familiar 
with the basics of Coxeter groups may safely skip this subsection.

Let $W$ be a finite reflection group in $\mathbf{GL}(n,\RR)$ and let $T$ be the family of reflections in $W$.
Denote by $\Phi$ the set of unit vectors orthogonal to the reflection hyperplanes, where each hyperplane gives rise to two vectors $\alpha$ and $-\alpha$.
This set of unit vectors $\Phi$ is called a \emph{root system} and each element $\alpha\in\Phi$ is called a \emph{root}.
Fixing an ordered basis $\{u_1,u_2,\cdots, u_n\}$ of $\RR^n$, 
we say that a root $\alpha$ is \emph{positive} if $\alpha=\sum c_i u_i$ and $c_k>0$, 
where $k$ is the smallest index $i$ for which $c_i\neq 0$.
Otherwise, a root is said to be \emph{negative}.
Clearly, the set $\Phi^+$ of positive roots and $\Phi^-$ of negative roots partition $\Phi$ and are of equal size.
Let $\Delta$ be a minimal subset of $\Phi^+$ such that each $\alpha\in\Phi^+$ is a linear combination 
of positive roots in $\Delta$ with non-negative coefficients.
Such a minimal subset always exists, since $\Phi^+$ itself already satisfies the condition. 
We call this $\Delta$ a \emph{simple system} and its elements are called \emph{simple roots}.
Given a simple root $\alpha$, the hyperplane orthogonal to $\alpha$ is called 
a \emph{simple reflection hyperplane}.
We denote by $s_\alpha$ the reflection through the simple reflection hyperplane orthogonal to $\alpha$
and refer to such reflections as \emph{simple reflections}.

The following theorem states some important facts about simple systems.
For a proof, we refer the reader to Sections 1.3 and 1.5 of \cite{H92}.
\begin{theorem}\label{thm:simple_system}
	A simple system $\Delta$ has the following properties:
	\begin{enumerate}[(i)]
		\item it is unique with respect to $\Phi^+$;
		\item it consists of linearly independent vectors;
		\item the set of all simple reflections generates $W$.
	\end{enumerate}		
\end{theorem}

For each positive root $\alpha$, let $H_\alpha$ be the reflection hyperplane orthogonal to $\alpha$.
Consider the collection $\mathcal{C}$ of connected components of $\RR^n\setminus\cup_{\alpha\in \Phi^+}H_\alpha$.
Each component $C$ in $\mathcal{C}$ is called an \emph{open chamber} of the reflection group $W$.
Denote by $\overbar{\C}=\{\overbar{C}:C\in\C\}$ the set of closures of open chambers, which we call \emph{closed chambers}.
Each open chamber consists of those vectors $v$ in $\RR^n$ with a certain
fixed sign for $\ang{v,\alpha}$ for each positive root $\alpha$.
Conversely, if we fix a sign for $\ang{v,\alpha}$ for each positive root $\alpha$, 
then, provided these choices are consistent, the collection of such vectors is an open chamber.
In particular, there is a unique open chamber consisting of all vectors $v$ satisfying $\ang{v,\alpha}>0$ for all positive roots $\alpha$,
since it is possible to have $\ang{v,\alpha}>0$ for all simple roots $\alpha$
and any such vector must have a positive inner product with every positive root. 
We call this chamber $C_0$ (or its closure $\overbar{C}_0$) the \emph{fundamental open (or closed) chamber}.
Note that the fundamental closed chamber $\overbar{C}_0$ is a cone given by 
	the intersection of closed half-spaces obtained from simple roots:
	\begin{align*}
		\overbar{C}_0=\bigcap_{\alpha\in\Delta}\{v\in\RR^n:\ang{v,\alpha}\geq 0\}.
	\end{align*}
	In other words, it is a closed cone surrounded by simple reflection hyperplanes,
	so a point $v\in \overbar{C}_0$ must be contained in either a simple reflection hyperplane
	or the open chamber $C_0$.

An important fact is that the action of $W$ on $\mathcal{C}$ is simply transitive.
We refer the reader to Sections 1.6 and 1.7 of \cite{H92} for a proof.
\begin{theorem}\label{thm:simple_trans}
	The action of $W$ on $\C$ is simply transitive.
	In particular, the identity is the only element in $W$ that fixes the fundamental open chamber.
\end{theorem}
It immediately follows that the action of $W$ on the set $\overbar{\C}$ of all closed chambers is also simply transitive.
For each reflection $t\in T$, let $H(t)$ be the hyperplane in $\mathbb{R}^n$ defining $t$ and let $D^+(t)$ be the open half-space $\{x\in\RR^n:\ang{x,\alpha}>0\}$, 
where $\alpha$ is the positive root orthogonal to $H(t)$.
Let $S\subset T$ be the set of all simple reflections and $I$ be a subset of $S$.
Define the cone $C(I)$ by
\begin{align*}
 C(I):= \left(\bigcap_{s\in I}H(s)\right)\cap
				\left(\bigcap_{s\in S\setminus I}D^+(s)\right).
\end{align*}
Observe that the fundamental open chamber $C_0$ is $C(\emptyset)$ and $C(I)$ and $C(J)$ are disjoint whenever $I\neq J$. 
Moreover, the fundamental closed chamber can be expressed as
\begin{align*}
	\overbar{C}_0=\bigcup_{I\subset S} C(I).
\end{align*}
The simple transitivity of the $W$-action may be extended as follows:
\begin{proposition}\label{prop:one_from_orbit}
The cone $C(I)$ is mapped into the fundamental closed chamber by $w\in W$ 
if and only if $w$ is in the subgroup of $W$ generated by $I$.
In particular, the stabiliser of $C(I)$ is precisely the subgroup generated by $I$.
\end{proposition}
Although we refer the reader to Theorem 3A8 in \cite{MS02} for more information, we remark that one direction of the proposition above is easy to see: if $w=s_1s_2\cdots s_k$ for $s_i\in I$,
then $w$ fixes each point in $C(I)$ since 
the subspace $\bigcap_{s\in I}H(s)$ is fixed under each $s_i$. Note, therefore, that if $C(I)$ is mapped into the fundamental closed chamber by an element of $W$, it is in fact mapped to itself.

\begin{example}\label{ex:S_n}
Let $W\cong S_n$ be the reflection group in $\mathbf{GL}(n,\RR)$ with reflection hyperplanes 
\begin{align*}
	H_{ij}:=\{(x_1,x_2,\cdots,x_n):x_i=x_j\}
\end{align*}
for $1\leq i<j\leq n$. Note that the reflection $t_{ij}$ through $H_{ij}$ is the map swapping $x_i$ and $x_j$.
With respect to the standard ordered basis $\{\mathbf{e}_1,\mathbf{e}_2,\cdots,\mathbf{e}_n\}$, the set of positive roots is 
\begin{align*}
	\Phi^+=\left\{\frac{1}{\sqrt{2}}(\ee_i-\ee_j):1\leq i<j\leq n\right\}.
\end{align*} 
The simple system $\Delta$ is then given by
\begin{align*}
	\Delta=\left\{\frac{1}{\sqrt{2}}(\ee_i-\ee_{i+1}):1\leq i\leq n-1\right\}
\end{align*}
and each permutation $\sigma\in S_n$ corresponds to an open chamber
\begin{align*}
	C_\sigma = \{(x_1,\cdots, x_n):x_{\sigma(1)}<x_{\sigma(2)}<\cdots<x_{\sigma(n)} \}.
\end{align*}
It is clear that $W$ is a simply transitive action on the set $\C=\{C_\sigma:\sigma\in S_n\}$ of all open chambers.
\end{example}

Each $w\in W$ can be expressed as $w=s_1s_2\cdots s_k$, where $s_1,s_2,\cdots,s_k\in S$.
Define the \emph{length} $\ell(w)$ of $w\in W$ to be the minimum length $k$ over all such expressions.
This purely combinatorial concept of length has a geometric interpretation,
which will play a crucial role in the next subsection.
\begin{theorem}\label{thm:geometric_length}
	Let $w\in W$ and $v\in\RR^n$ be a point in the fundamental open chamber.
	For each positive root $\alpha$ and the reflection $t_{\alpha}\in T$ induced by $\alpha$, 
	$\ell(t_\alpha w)>\ell(w)$ holds if and only if $\ang{w(v),\alpha}$ is positive. 
\end{theorem}

The theorem above says that 
$t_\alpha w$ has a greater length than $w$ if and only if
$w(v)$ lies on the `positive side' of $H_\alpha$, i.e., the same side as the positive root $\alpha$, whereas its image $t_\alpha w(v)$ under the reflection $t_\alpha$ lies on the other side.
Observe now that a root $\alpha$ is positive if and only if $\ang{v,\alpha}>0$
whenever $v$ is chosen from the fundamental open chamber.
Thus, the `positive side' with respect to the hyperplane $H_\alpha$ is exactly 
the component of $\RR^n\setminus H_\alpha$ containing the fundamental open chamber.
To summarise, Theorem \ref{thm:geometric_length} implies that the length of $t_\alpha w$ becomes greater than $w$ if and only if by applying $t_\alpha$ we move $w(v)$ 
from the same side of $H_\alpha$ as the fundamental open chamber
to the opposite side.
Theorem \ref{thm:geometric_length} is a key fact in the theory of Coxeter groups,
so proofs can be found in many places: 
for instance, Sections 1.6 and 1.7 in \cite{H92} or Proposition 4.4.6 in \cite{BB05}.

\subsection{The Euclidean embedding} \label{sec:Euclidean}

We begin with a motivating example:
\begin{example}\label{ex:K_4}
Let $H$ be the 1-subdivision of $K_4$. This graph is \emph{biregular}, that is, it is a bipartite graph such that all vertices
on the same side of the bipartition have the same degree. In this case, all vertices on one side have degree three, while 
all vertices on the other side have degree two.
One may check that the cut involution group $W_H$ of $H$ is isomorphic to the symmetric group $S_4$.
Note also that $S_4$ may be represented as the reflection group of the tetrahedron.
There is a natural way to embed $H$ into $\mathbb{R}^4$ that uses these observations:
consider the map $\eta:V(H)\rightarrow\mathbb{R}^4$ which sends the vertices $v_1,v_2,v_3$, and $v_4$ of degree 3 to $\mathbf{e}_1,\mathbf{e}_2,\mathbf{e}_3$, and $\mathbf{e}_4$, respectively, 
and the unique common neighbour $u_{ij}$ of $v_i$ and $v_j$, $i\neq j$, to the midpoint $\frac{1}{2}(\mathbf{e}_i+\mathbf{e}_j)$ of their images under $\eta$.
Then every reflection in $S_4$ becomes a cut involution of $H$ and, 
conversely, every cut involution is represented by a reflection.
\end{example}

This example shows that some graphs may be embedded in Euclidean space so that cut involutions are represented by genuine reflections.
Our aim in this subsection will be to confirm the existence of a similar embedding for every reflection graph.

We will continue to use notation from Section \ref{sec:Prelim}. That is, we let $W\subset\mathbf{GL}(n,\RR)$ be a finite reflection group, $T$ be the family of reflections in $W$, and $S\subset T$ be the set of simple reflections with respect to some fixed ordered basis of $\RR^n$.
Fix also subsets $S_1$ and $S_2$ of $S$ and let $W_1$ and $W_2$ be the subgroups of $W$ generated by $S_1$ and $S_2$, respectively. For brevity in what follows, we will write $H$ for the $(S_1, S_2; S, W)$-reflection graph. Recall that this is the bipartite graph between the (left-)cosets of $W_1$ and $W_2$,
where $wW_1$ and $wW_2$ are adjacent for every $w\in W$.

In Example~\ref{ex:K_4}, we mapped the vertices of $H$ to points in $\RR^n$ so as to visualise the graph more easily.
However, to avoid unnecessary ambiguity in the general case,
we will instead construct a map $\eta$ that sends vertices of $H$ to cones in $\RR^n$.
Indeed, let $\eta$ be the map from $V(H)$ to the family $\{wC(I):w\in W, I\subset S\}$ of cones in $\RR^n$ such that,
for $i=1,2$ and each $w\in W$,
\begin{align*}
	\eta(wW_i):=wC(S_i).
\end{align*}
In particular, $\eta(W_1)=C(S_1)$ and $\eta(W_2)=C(S_2)$, i.e., they are mapped into subcones of the closed fundamental chamber.
Since $uW_i=wW_i$ if and only if $w^{-1}u\in W_i$ and the cone $C(S_i)$ is fixed by each $g\in W_i$, we have
\begin{align*}
	\eta(uW_i)=uC(S_i)=w(w^{-1}u)C(S_i)=wC(S_i)=\eta(wW_i),
\end{align*}
and hence the map $\eta$ is well-defined. We call the map $\eta$ the \emph{Euclidean embedding} of $H$.

\begin{example}\label{ex:S_4}
It will be instructive to revisit the case where $H$ is the 1-subdivision of $K_4$.
Following the notation in Example~\ref{ex:S_n} with $n=4$, let $W\cong S_4$ be the reflection group in $\mathbf{GL}(4,\RR)$ 
with the reflection $t_{ij}$ swapping $x_i$ and $x_j$, $1\leq i<j\leq 4$.
Let $S_1=\{t_{23},t_{34}\}$ and $S_2=\{t_{12},t_{34}\}$.
Then the $(S_1,S_2;S,W)$-reflection graph $H$ is isomorphic to the $1$-subdivision of $K_4$.
The Euclidean embedding $\eta$ is then given by
\begin{align*}
	\eta(W_1)=\{(x_1,x_2,x_3,x_4):x_1> x_2=x_3=x_4\},\\
	\eta(W_2)=\{(x_1,x_2,x_3,x_4):x_1=x_2> x_3= x_4\},
\end{align*}
and the other values are given by permuting some coordinates of $\eta(W_1)$ and $\eta(W_2)$.
Observe that $W_1$ and $W_{2}$ are the only pair of adjacent vertices 
such that both of their images under $\eta$ are subsets of the fundamental closed chamber $\{(x_1,x_2,x_3,x_4):x_1\geq x_2\geq x_3\geq x_4\}$.
One may also check that the points used in Example~\ref{ex:K_4} to embed vertices are `typical' points from the corresponding cones.
\end{example}

The key properties of the Euclidean embedding are spelled out in the following lemma. 

\begin{proposition}\label{prop:geometric_structure}
Let $H$ be the $(S_1,S_2;S,W)$-reflection graph.
Then the Euclidean embedding $\eta$ of $H$ has the following properties:
	\begin{enumerate}[(i)]
		\item for two vertices $a$ and $a'$ on the same side of the bipartition of $H$, $\eta(a)=\eta(a')$ if and only if $a=a'$;
		\item every closed chamber $\overbar{C}$ contains both $\eta(a)$ and $\eta(b)$ as subsets for exactly one edge $ab\in E(H)$;
		\item there exists a closed chamber $\overbar{C}$ containing both $\eta(a)$ and $\eta(b)$ whenever $ab$ is an edge in $H$.
	\end{enumerate}
\end{proposition}
\begin{proof}
	\begin{enumerate}[(i)]
		\item 
		Suppose that $a=wW_1$, $a'=uW_1$ and $\eta(a)=\eta(a')$.
		By the definition of $\eta$, we have $\eta(a)=wC(S_1)$ and $\eta(a')=uC(S_1)$.
		Since $\eta(a)=\eta(a')$, it follows that $uw^{-1}$ fixes $C(S_1)$.
		By Proposition \ref{prop:one_from_orbit},
		this implies that $uw^{-1}\in W_1$, and hence $a=a'$.
		
		\item It follows from the construction of $\eta$ that the fundamental closed chamber $\overbar{C}_0$ contains $\eta(W_1)$ and $\eta(W_2)$.
		Again by Proposition \ref{prop:one_from_orbit}, $\eta(W_1)$ is mapped to a cone in $\overbar{C}_0$ by $w\in W$ if and only if $w\in W_1$ 
		and hence every point in $\eta(W_1)$ is mapped to itself.
		Since $\eta$ is one-to-one on each side of the bipartition by (i), 
		$W_1$ and $W_2$ are the only vertices of $H$ that are mapped to subsets of $\overbar{C}_0$.
		Moreover, each closed chamber $w\overbar{C}_0$ contains exactly one pair of cones $(\eta(wW_1),\eta(wW_2))$ that are images of an adjacent pair 
		$(wW_1, wW_2)$, since 
		the action of $W$ is simply transitive on $\C$ and the fundamental closed chamber contains $\eta(W_1)$ and $\eta(W_2)$.
		
		\item Since $ab\in E(H)$, there exists $w \in W$ such that $a=wW_1$ and $b=wW_2$.
		Then $\eta(a)$ and $\eta(b)$ must be contained in $w\overbar{C}_0$, where $\overbar{C}_0$ is the fundamental closed chamber.\qedhere
	\end{enumerate}
\end{proof}

The next corollary highlights the fact that we have identified a class of cut involutions in $H$, namely, those corresponding to genuine reflections in the Euclidean embedding.

\begin{corollary}
	Let $H$ be the $(S_1,S_2;S,W)$-reflection graph. Then every reflection $t\in T$ is a cut involution of $H$.
\end{corollary}
\begin{proof}
	It is straightforward to check that every $t\in T$ acts as an involutary graph automorphism on $H$. 
	By part (iii) of the proposition above, every edge in $H$ is mapped into a closed chamber by $\eta$, and hence there is no edge crossing a reflection hyperplane.
\end{proof}

One may wonder if the converse of the corollary above also holds, that is, whether
every cut involution of a reflection graph $H$ becomes a reflection under the Euclidean embedding.
However, this is not true in general, as evidenced by the following example:
\begin{example}\label{ex:K1,4}
Suppose that $W\cong S_4$, $S=\{t_{12},t_{23},t_{34}\}$, and $S_1=\{t_{23},t_{34}\}$, as in Example \ref{ex:S_4},
but take $S_2=S$ instead.
Then the $(S_1,S_2;S,W)$-reflection graph is isomorphic to the star $K_{1,4}$. However, it is not possible to render all 
cut involutions so that they correspond to reflections. This is because $K_{1,4}$ has a cut involution which is a product of two other cut involutions and, therefore, the determinants of all three of the corresponding linear transformations cannot be $-1$.
\end{example}

We have shown that every reflection $t \in T$ induces a cut involution on $H$, mapping the vertex $wW_i$ to the vertex $twW_i$. Since $T$ generates $W$, the graph $H$ is edge-transitive under the action of these cut involutions. Thus, in order to apply Theorem~\ref{cor:colouring_graphs}, it remains to show that there exists a percolating sequence. This will be the topic of the next subsection.

\subsection{Proof of Theorem \ref{thm:parabolic}}
\begin{figure}
    \centering
    \includegraphics[width=0.4\textwidth]{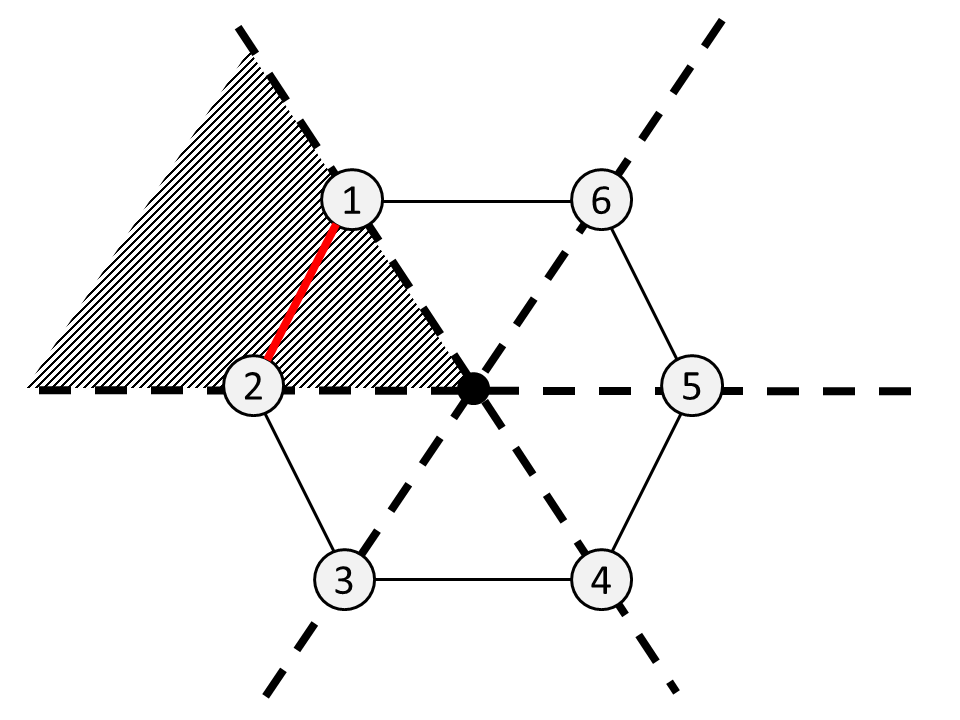}
    \caption{The Euclidean embedding of $C_6$.}\label{fig:C6}
\end{figure}

Proposition \ref{prop:geometric_structure} allows us to lift
the notion of folding sequences to a group theoretic setting.
Here is an illustrative example:
\begin{example}
Recalling the example of $C_6$ on vertex set $\{1,2,\cdots,6\}$ discussed in Section \ref{sec:C6},
suppose that we have oriented each cut involution in such a way that
the component containing the edge $\{1,2\}$ is on the left.
Then the folding sequence starting from the single edge set $J_0:=\{\{1,2\}\}$ with $J_{i+1}=J_i^+(\phi_i)$ at the $i$-th step,
where the $\phi_i$ are chosen as in Section \ref{sec:C6}, 
is a percolating sequence.
Now consider $C_6$ as the $(S_1,S_2;S,W)$-reflection graph, where $W=I_2(3)\subset\mathbf{GL}(2,\RR)$ is the dihedral group,
$S=\{s_1,s_2\}$ is a generating set of simple reflections of $I_2(3)$ with angle $\pi/3$ between the reflection lines defining $s_1$ and $s_2$,
$S_1=\{s_1\}$, and $S_2=\{s_2\}$.
Figure \ref{fig:C6} shows what $\eta$ looks like if we label $W_1$ and $W_2$ with $1$ and $2$, respectively: 
we have chosen a basis of $\RR^2$ in such a way that the fundamental open chamber $C_0$, represented by the grey region, is the unique chamber
whose closure contains $\eta(1)$ and $\eta(2)$ as subsets
and each vertex is mapped to a ray from the origin as marked. 
The folding sequence described above naturally gives a folding sequence on chambers, allowing us to spread the grey colour on 
the fundamental chamber across all chambers. For instance, we may spread the colour on the grey chamber in Figure \ref{fig:C6}
to the chamber containing the edge $\{2,3\}$ by applying the reflection through the horizontal line through $2$ and $5$. In this subsection,
we will approach the general question in the opposite direction, using the algebraic tools built up over the previous subsections to show how to spread the grey colour across all chambers, and then using this to form a folding sequence for the original graph.
\end{example}

Let $\eta$ be the Euclidean embedding of the $(S_1,S_2;S,W)$-reflection graph $H$.
For each reflection $t\in T$ and the associated positive root $\alpha$, let $t^+$ and $t^-$ be 
the left and right-folding maps defined by
\begin{align*}
	t^+(x):=
	\begin{cases}
		t(x) &\text{ if }\ang{x,\alpha}< 0\\
		x &\text{ if }\ang{x,\alpha}\geq 0
	\end{cases}\text{ and ~}
	t^-(x):=
	\begin{cases}
		x &\text{ if }\ang{x,\alpha}\leq 0\\
		t(x) &\text{ if }\ang{x,\alpha}> 0,
	\end{cases}
\end{align*}
respectively.
Given a family $\J$ of open chambers, we set
\begin{align*}
	\J^+(t):=\{C\in \C:t^+(C)\in \J\}\text{ and }
	\J^-(t):=\{C\in\C:t^-(C)\in \J\}.
\end{align*}
Again, as for the graph case, we say that a sequence $\J_0,\J_1,\cdots$ of families of open chambers
is a \emph{folding sequence} in $\C$ if for each $i\geq 0$ there exists $t\in T$ such that 
\begin{align*}
	\J_{i+1}=\J_i^+(t) \text{ or }\J_{i+1}=\J_i^-(t).
\end{align*}
Similarly, we say that a finite folding sequence $\J_0,\J_1,\cdots,\J_N$ is a \emph{percolating sequence} 
if it starts with a single chamber $\J_0$ and ends with the set of all open chambers $\J_N=\C$.

We claim that the existence of a percolating sequence in $\C$ implies 
the existence of a percolating sequence in the corresponding graph $H$.
Suppose that there exists a folding sequence $\J_0,\cdots,\J_N=\C$.
By the transitivity of the $W$-action on $\C$, 
we may assume that $\J_0$ is the fundamental open chamber $C_0$.
At the $i$-th step, we say that an open chamber $C$ or its closure $\overbar{C}$ is coloured if $C\in\J_i$.
As chambers are coloured, we project the colouring down to edges of $H$
by regarding an edge $ab$ as coloured if both $a$ and $b$ are embedded in the 
closure of a coloured open chamber.
This projection respects the folding operations:
if the colour on a chamber $C$ spreads to $t(C)$ under the reflection $t$, 
then the colour on the edge $ab$ with $\eta(a),\eta(b)\subset \overbar{C}$
spreads to the edge $t(a)t(b)$.

Formally, 
let $J_\J\subset E(H)$ be the set of edges $ab$ such that $\eta(a),\eta(b)\subset  \overbar{C}$ for some $C\in\J$.
Recall that each reflection $t\in T$ acts in two different ways:
it acts on the set of chambers $\C$ group theoretically and on the vertex set $V(H)$ consisting of all cosets of $W_1$ and $W_2$ as a cut involution.
To distinguish these two actions, denote by $\phi_t$ the cut involution of $H$ that corresponds to $t\in T$.   
We will show a correspondence
\begin{align}\label{eqn:correspondence}
	J_{\J^+(t)}=J_\J^+(\phi_t)\text{ and }J_{\J^-(t)}=J_\J^-(\phi_t)
\end{align}
between edge sets and collections of open chambers,
provided that each cut involution $\phi_t$ is oriented so as to guarantee the consistency of signs. 
By this correspondence, a percolating sequence $\{C_0\}=\J_0,\cdots,\J_N=\C$ induces 
a folding sequence $J_0,\cdots,J_N$ with $J_i=J_{\J_i}$.
Recall that, by Proposition \ref{prop:geometric_structure} (ii),
 $(W_1,W_2)$ is the unique pair of vertices both of whose images under $\eta$ are contained in the fundamental closed chamber.
Thus, $J_0$ must consist of the single edge $\{W_1,W_2\}$. Since we also have $J_N=E(H)$ by Proposition \ref{prop:geometric_structure} (iii), 
$J_0,\cdots,J_N$ is the desired percolating sequence.

It remains to prove the correspondence \eqref{eqn:correspondence}.
We orient each cut involution $\phi_t$ in such a way that
\begin{align*}
	L_{\phi_t}:=\{x\in V(H): \eta(x)\subset D^+(t))\}\text{ and }
	R_{\phi_t}:=\{x\in V(H): \eta(x)\subset D^-(t)\}.
\end{align*}
For brevity, we write $L_t:=L_{\phi_t}$ and $R_t:=R_{\phi_t}$.
Observe that $L_{t}$ and $R_{t}$ are disjoint sets which are mapped to each other by $t$
and the remaining vertices $x\in V(H)\setminus(L_{t}\cup R_{t})$ 
with $\eta(x)\subset H(t)$ are fixed points of $t$.
With this orientation, the left-folding map $\phi_t^+$ sends $x$ to $y$ if and only if $t^+$ maps $\eta(x)$ onto $\eta(y)$.
Therefore, we have
\begin{align*}
	J_\J^+(\phi_t)&=\{ab\in E(H):\phi_t^+(a)\phi_t^+(b)\in E(J_\J)\}\\
			&=\{ab\in E(H):t^+(\eta(a)\cup\eta(b))\subset \overbar{C}\text{ for some }C\in\J\}\\
			&=J_{\J^+(t)}
\end{align*}
and, similarly, $J_\J^-(\phi_t)=J_{\J^-(t)}$.

\medskip

We will now prove the existence of percolating sequences in $\C$ by reducing to a purely group-theoretic framework.
Let $S$ be the set of simple reflections in $W$ given by the fixed simple system $\Delta$
and let $p$ be a point in the fundamental open chamber. 
There is a natural one-to-one correspondence between open chambers and elements of $W$ that maps $C\in\C$ to the unique element $w\in W$ such that $w(p)\in C$. We will write $C=wC_0$ if $w(p)\in C$.
For example, the identity element of $W$ corresponds to the fundamental open chamber $C_0$.
For each reflection $t$, we may define the left and right-folding maps $t_+$ and $t_-$ on $W$ so that
\begin{align*}
	(t_+ w)C_0=t^+(wC_0)\text{ and }(t_- w)C_0=t^-(wC_0)
\end{align*}
hold.
Recall that $t^+(C)=t(C)$ if and only if $\ang{x,\alpha}< 0$ for all $x\in C$,
where $\alpha$ is the positive root that induces $t$.
In other words, $t_+w=tw$ if and only if for every $x\in wC_0$ we have $\ang{x,\alpha}<0$.
By Theorem \ref{thm:geometric_length}, this is equivalent to the purely algebraic property $\ell(tw)<\ell(w)$.
Thus, one may check that
\begin{align}\label{eqn:half_fold_group}
	t_+ w=
	\begin{cases}
		tw \text{~ if }\ell(tw)<\ell(w)\\
		w \text{~~ if }\ell(tw)\geq\ell(w)
	\end{cases}\text{ and ~}
	t_-w=
	\begin{cases}
		w \text{~~ if }\ell(tw)\leq\ell(w)\\
		tw \text{~ if }\ell(tw)>\ell(w)
	\end{cases}
\end{align}
hold.\footnote{In fact, it is impossible to have $\ell(tw)=\ell(w)$. 
This follows from the simple fact that $w\mapsto \det(w)=(-1)^{\ell(w)}$ is a group homomorphism from $W$ to the multiplicative group $\{\pm1\}$.}
For a subset $K\subset W$, let
\begin{align*}
	K^+(t):=\{w\in W:t_+w\in K\}\text{ and }
	K^-(t):=\{w\in W:t_-w\in K\}.
\end{align*}
For $\J\subset\C$, write $K_\J$ for the subset $\{w\in W:wC_0\in\J\}$ of $W$. 
We have a natural correspondence $K_{\J^+(t)}=K_\J^+(t)$ and $K_{\J^-(t)}=K_\J^-(t)$.
Again, we say that a sequence of subsets $K_0,K_1,\cdots$ is 
a folding sequence in $W$ if and only if
$K_{i+1}$ equals $K_i^+(t)$ or $K_i^-(t)$ for each $i\geq 0$.
Hence, it remains to find a folding sequence $\{e\}=K_0,K_1,\cdots,K_N=W$. The existence
of such a sequence follows from the next theorem. In fact, it shows something stronger:
we may always choose $K^+(s)$ with $s$ a simple reflection to update $K$.
\begin{theorem}\label{thm:colouring_chambers}
	Let $W$ be a finite reflection group and $e$ the identity of $W$.
	Then there exists a finite folding sequence $\{e\}=K_0,K_1,\cdots,K_N=W$ such that
	for each $i=0,1,\cdots,N-1$ there is a simple reflection $s_i$ for which $K_{i+1}=K_i^+(s_i)$. 
\end{theorem}
Our plan for proving Theorem~\ref{thm:colouring_chambers} is to use induction on the length function $\ell$.
We say that a subset $U$ of $W$ is a \emph{stack} 
if whenever $w\in U,s\in S$, and $\ell(sw)<\ell(w)$, $sw$ is also in $U$.\footnote{That is, 
a stack is a subset of $W$ closed under downward inclusion with respect to the Bruhat (or strong) order, which is the poset structure induced by length.
We have not defined the Bruhat order, but we use it implicitly.}
The following lemma is the key to proving Theorem \ref{thm:colouring_chambers}.
\begin{lemma}\label{lem:stack}
 Let $U$ be a stack and $s\in S$ be a simple reflection. 
 Then $U$ is contained in $U^+(s)$.
\end{lemma}
\begin{proof}
	Let $u\in U$ be arbitrary. 
	If $\ell(su)<\ell(u)$, the fact that $U$ is a stack immediately implies  that $su\in U$.
	By \eqref{eqn:half_fold_group}, we then have $s_+u=su$ and hence $s_+u\in U$,
	which means $u\in U^+(s)$.
	Otherwise, $s_+u=u\in U$, which again implies $u\in U^+(s)$.
\end{proof}

\begin{proof}[Proof of Theorem \ref{thm:colouring_chambers}]
	Let $S=\{s_1,\cdots,s_k\}$ be the set of simple reflections, which is a generating set for $W$.
	We will prove by induction on $L$ that there is a $K_i$ which contains all $w\in W$ with $\ell(w)\leq L$.
	The only element with length $0$ is the identity, so $K_0 = \{e\}$ satisfies the condition.
	Suppose now that $U_L$ is the set of all $w\in W$ with $\ell(w)\leq L$ and $K_{i_L}$ contains $U_L$.
	Let $U_{L+1,0}=W_{L+1,0}=U_L$ and define, for $i,j=1,2,\cdots,k$,
	\begin{align*}
		U_{L+1,i}:=\{w\in W:s_{i}w\in U_L \}\text{ and }W_{L+1,j}:=\bigcup_{i=0}^j U_{L+1,i}.
	\end{align*}
	Note that $W_{L+1,j}$ is always a stack, because elements of smaller length must be in $U_{L}$. 
	Therefore, by Lemma \ref{lem:stack}, $W_{L+1,j}$ is a subset of $W_{L+1,j}^+(s_{j+1})$ for each $j=0,1,\cdots,k-1$. Moreover, $U_{L+1, j+1}$ is a subset 
	of $W_{L+1,j}^+(s_{j+1})$ for each $j=0,1,\cdots,k-1$. To see this, note that $w \in U_{L+1, j+1}$ if and only if $s_{j+1} w \in U_L$. If $w \in U_L$, $w$ is 
	already contained in $W_{L+1,j}$. We may therefore assume that $w \not\in U_L$ and so $\ell(s_{j+1} w) < \ell(w)$. This in turn shows that 
	$(s_{j+1})_+ w = s_{j+1} w \in W_{L+1,j}$, which implies that $w \in W_{L+1,j}^+(s_{j+1})$.
	Putting everything together, we see that $W_{L+1,j+1}$ is a subset of $W_{L+1,j}^+(s_{j+1})$ for each $j=0,1,\cdots,k-1$.
	Since $W_{L+1,k}$ is the set $\{w\in W:\ell(w)\leq L +1\}$, we may therefore take $i_{L+1} = i_L + k$, completing the induction.
\end{proof}

Since we have now shown the existence of a percolating sequence in $H$, this also completes the proof of Theorem~\ref{thm:parabolic}. Theorem~\ref{thm:polytope_norming}, concerning incidence graphs of regular polytopes, follows as a simple application.

\begin{proof}[Proof of Theorem \ref{thm:polytope_norming}]
	Let $\cP$ be a (realised) regular polytope with symmetry group $W$
	and let $\Psi$ be a fixed flag in $\cP$, that is, a maximal chain of faces.
	It is a folklore fact (see, for example, Theorem 3D7 in \cite{MS02}) that we may choose a basis of $\RR^n$
	so that the set of simple reflections $S$ can be enumerated as $\{s_1,\cdots,s_n\}$ in such a way that 
	the $k$-face $F_k$ in $\Psi$ is fixed under all reflections in $S$ but $s_{k+1}$,
	i.e., $s_{k+1}$ is the only active mirror for $F_k$ in $S$.
	Since $\cP$ is regular, there is a one-to-one correspondence between 
	the cosets of the parabolic subgroup generated by $S\setminus\{s_{k+1}\}$ and the $k$-faces of $\cP$.
	Hence, it follows that the $(S_1,S_2;S,W)$-reflection graph with $S_1=S\setminus \{s_{k+1}\}$ and $S_2=S\setminus \{s_{r+1}\}$
	is isomorphic to the $(k,r)$-incidence graph of $\cP$.
\end{proof}

There are many more examples than those coming from regular polytopes.
For instance, the vertex-edge incidence graphs of many quasiregular, i.e., vertex and edge-transitive, polytopes can also be written as $(S_1,S_2;S,W)$-reflection graphs.
Moreover, it is possible to obtain interesting examples by considering the orbits of faces of the same rank:

\begin{example}\label{ex:hypercube}
Let $D_n$ be the $n$-demicube group, i.e., the order $2$ subgroup of the hypercube group obtained by removing the `perpendicular' mirrors $x_i=0$.
The action of $D_n$ on the set of vertices of an $n$-dimensional hypercube
is not transitive, but there exist two orbits,
each of which corresponds to a colour class of the hypercube when considered as a bipartite graph. 
There is a generating set of simple reflections $S=\{s_1,s_2,\cdots,s_n\}$ for $D_n$ such that
the order of $s_1s_3$ and $s_is_{i+1}$ for $i\geq 2$ is $3$, while all other $s_is_j$, $i\neq j$, are of order $2$.
Then one may check that the $n$-dimensional hypercube is isomorphic to the
$(S_1,S_2;S,D_n)$-reflection graph with $S_1=S\setminus\{s_1\}$ and $S_2=S\setminus \{s_2\}$. 
It follows that the hypercube is weakly norming, a fact first proved in \cite{H10}.
\end{example}
By Theorem \ref{thm:polytope_norming}, we already know that the vertex-edge incidence graph of an octahedron,
i.e., the 1-subdivision of an octahedron, is weakly norming, but we may also prove that it is norming through a judicious choice of reflection group.
\begin{example}\label{ex:norming_D3}
	Let $D_3$ be the 3-dimensional demicube group and $S$ the same generating set described in Example \ref{ex:hypercube}.
	Take $S_1=\{s_1,s_2\}$ and $S_2=\{s_3\}$. 
	Then $(S_1,S_2;S,D_3)$-reflection graph is isomorphic to the 1-subdivision of an octahedron,
	as the latter is the face-edge incidence graph of the cube.
	Therefore, the $1$-subdivision of the octahedron is norming.
\end{example}
Let us mention another example of a norming graph.
Following \cite{CKLL15}, 
we say that a graph is a $K_{2,t}$-replacement of $H$ if each edge of $H$ is replaced with 
a copy of $K_{2,t}$ by identifying the two vertices of the edge with the two vertices on the smaller side of $K_{2,t}$.
\begin{example}\label{ex:norming_B3}
	Let $B_3$ be the cube group, that is, the symmetry group of the cube.
	Then there is a generating set $S=\{s_1,s_2,s_3\}$ such that $s_1s_2$, $s_2s_3$, and $s_1s_3$ are of orders $4$, $3$, and $2$, respectively.
	Let $S_1=\{s_1,s_2\}$ and $S_2=\{s_3\}$.
	Then one may check that the $(S_1,S_2;S,B_3)$-reflection graph is isomorphic to the $K_{2,2}$-replacement of the octahedron graph,
	and thus is norming.
\end{example} 
Finally, we remark that one may presumably use the exceptional reflection groups $E_6,E_7$, and $E_8$ to build some more exotic (weakly) norming graphs, though we have not pursued this further.

\section{Generalisations and applications} \label{sec:genapp}
\subsection{Hypergraph norms}

As noted in Hatami's PhD thesis~\cite{H09}, the concepts of norming and weakly norming graphs generalise in the obvious way to hypergraphs, 
with Gowers' octahedral norms \cite{G06,G07} serving as standard examples.
In this short subsection, we discuss the appropriate generalisation of Theorem \ref{thm:parabolic}.

To obtain a suitable generalisation of Theorem \ref{thm:parabolic} to $k$-uniform hypergraphs, or $k$-graphs for short, with the same proof strategy, 
we should first define cut involutions for hypergraphs.
One naive way might be to say that an involutory automorphism is a cut involution 
if the fixed point set is again a vertex cut, i.e., deleting it makes the hypergraph disconnected.
However, this is not a good choice.
For example, consider the 3-graph $H$ on vertex set $\{a,b,c,x,y,z\}$ with three edges $xbc, yca,$ and $zab$.
Let $\phi$ be the involutory automorphism of $H$ that fixes $x$ and $a$ and maps $b$ and $z$ to $c$ and $y$, respectively.
Deleting the fixed vertex set consisting of $x$ and $a$ obviously makes $H$ disconnected,
but we cannot use $x$ and $a$ as `pivots' to apply the Cauchy--Schwarz inequality, since there is an edge crossing the cut.
Furthermore, if $H$ were a (weakly) norming hypergraph then it would also satisfy the hypergraph generalisation of Sidorenko's conjecture, 
but $H$ was shown to be a counterexample to this conjecture in \cite{Sid92}. 

To define cut involutions so as to avoid the difficulties discussed in the previous example, we replace each edge in a $k$-graph $H$ by a clique of size $k$ and write $\tilde{H}$ for the resulting graph on $V(H)$.
Note that every automorphism $\phi$ of $H$ is also an automorphism of $\tilde{H}$.
We say that an automorphism $\phi$ of $H$ is a cut involution of $H$ if it is a cut involution of $\tilde{H}$.
With this definition, all of our arguments generalise without difficulty.

Given a finite reflection group $W$, the set of all simple reflections $S$, and subsets $S_1,S_2,\cdots,S_k$ of $S$,
let the $(S_1,\cdots,S_k;S,W)$-\emph{reflection hypergraph} be the $k$-partite $k$-graph 
whose parts are the cosets of the subgroup $W_i$ generated by $S_i$ for each $i=1,\cdots,k$, with an edge for every $k$-tuple of the form
$(wW_1,wW_2,\cdots,wW_k)$ with $w\in W$.
We say that a $k$-graph is a \emph{reflection hypergraph} if it is isomorphic 
to the $(S_1,\cdots,S_k;S,W)$-reflection hypergraph for a suitable choice of parameters.
We have the following hypergraph generalisation of Theorem \ref{thm:parabolic}:

\begin{theorem}\label{thm:hypergraph}
	A $k$-partite $k$-graph $H$ is weakly norming whenever it is a reflection hypergraph.
	Moreover, if there is no mirror containing a hyperedge, i.e., $\bigcap_{i=1}^k S_i=\emptyset$, then $H$ is norming.
\end{theorem}

With this framework, we may easily recover Gowers' octahedral norms~\cite{G06,G07}.

\begin{example}\label{ex:Gowers_norms}
	Let $W=W_1\times W_2\times\cdots\times W_k$, where each $W_i$ is a reflection group generated by a single reflection $s_i$.
	Then the set $S$ of simple reflections consists of $k$ orthogonal reflections.
	If we let $S_i=S\setminus\{s_i\}$ 
	and $H_k$ be the $(S_1,\cdots,S_k;S,W)$-reflection hypergraph,
	then Theorem \ref{thm:hypergraph} implies that $H_k$, consisting of the $(k-1)$-faces of the $k$-dimensional octahedron, is a norming $k$-graph.
\end{example}

Another example of a hypergraph norm comes from the work on weak quasirandomness in \cite{CHPS12}.

\begin{example}\label{ex:vertex-quasirandom}
	In \cite{CHPS12}, a $k$-graph $M_k$ is constructed recursively as follows.
	Given a $k$-partite $k$-graph $M$ on $A_1\cup\cdots\cup A_k$,
	we write $\db_i(M)$ for the $k$-graph obtained by gluing two vertex-disjoint copies of $M$ 
         so that the corresponding vertices in each copy of $A_i$ are identified.
	Starting from the graph $M_0$ with a single edge,
	define
	\begin{align*}
		M_k:=\db_k(\db_{k-1}(\cdots \db_1(M_0)\cdots)).
	\end{align*}
	Letting $S_i=\{s_i\}$ and using the same reflection group $W$ as in Example \ref{ex:Gowers_norms},
	Theorem \ref{thm:hypergraph} implies that the $(S_1,\cdots,S_k;S,W)$-reflection hypergraph is norming. 
	Since this hypergraph is isomorphic to $M_k$, we see that $M_k$ defines a semi-norm.
\end{example}

The examples above are very natural: both of them use the simplest reflection group with $k$ generators and the resulting Cauchy--Schwarz trees are very symmetric and easy to analyse.
Our new framework gives 
a much larger class of hypergraph norms, each of which defines a certain notion of quasirandomness.  
We will discuss the relations between these notions in the next subsection.

\subsection{Domination between norms} 

In this subsection, we will be interested in relations between (hyper)graph norms. 
To begin, we will study the question of determining whether
the (absolute) $H$-norm dominates another (absolute) $J$-norm.
That is, for any bounded measurable function $f:[0,1]^2 \rightarrow\RR$, we would like to know if
\begin{align} \label{eqn:mono}
	\norm{f}_J\leq\norm{f}_H ~\left(\text{or }\norm{f}_{r(J)}\leq\norm{f}_{r(H)}\right)
\end{align}
holds.
This question remains valid even if $H$ and $J$ are not (weakly) norming,
but weakly norming graphs can be regarded as local maxima for such comparisons.
To be more precise, whenever $H$ is (weakly) norming and $J$ is a subgraph of $H$,
the (absolute) $H$-norm dominates the (absolute) $J$-norm.
To see this, let $f_{\chi(e)}=1$ for all $e\notin E(J)$ in inequality \eqref{eqn:rainbow}. 

An immediate consequence of this inequality is Sidorenko's conjecture \eqref{eqn:Sidorenko},
which essentially states that the absolute $H$-norm dominates the single-edge norm whenever $H$ is a bipartite graph.
Hence, if $H$ is weakly norming, Sidorenko's conjecture holds for $H$. Similar arguments show that if $H$ is weakly norming and contains a cycle, then it also satisfies the so-called forcing conjecture, a central problem in the study of quasirandom graphs. We refer the interested reader to~\cite{CFS10} for further information about this conjecture and its relationship with graph norms.

We remark that one (hyper)graph norm may dominate another  
even when the first (hyper)graph does not contain the second. For example, if $m \geq n$, then
\begin{align*}
 \norm{f}_{C_{2m}}\leq\norm{f}_{C_{2n}}.
\end{align*}
Note that the cut involution groups of different even cycles are always non-isomorphic.
In general, it is hard to compare (hyper)graph norms coming from non-isomorphic reflection groups,
but we can say something if the two (hyper)graphs have isomorphic cut involution groups of a certain type.

\begin{proposition}\label{prop:domination}
	Let $W=W_1\times W_2\times\cdots\times W_n$, where each $W_i$ is a reflection group generated by a single reflection $s_i$, 
	and let $S=\{s_1,\cdots,s_n\}$ be the set of simple reflections in $W$.
	Suppose $S_1,\cdots,S_k$ and $S_1',\cdots,S_k'$ are subsets of $S$ 
	such that $S_i'\subset S_i$ for each $i=1,2,\cdots,k$.
	Let $H$ and $H'$ be the $(S_1,\cdots,S_k;S,W)$-reflection hypergraph and the $(S_1',\cdots,S_k';S,W)$-reflection hypergraph, respectively.
	Then both $H$ and $H'$ are norming and the $H$-norm dominates the $H'$-norm.
\end{proposition}
Our proof is quite heavy on notation, so we postpone it until the end of this subsection and instead give an example 
that conveys the rough idea of the proof.
\begin{example}\label{ex:domination}
Consider the special case of Proposition \ref{prop:domination} 
where $n=3$, $S_i=S\setminus\{s_i\}$ for $i=1,2,3$, $S_2'=S_2$, $S_3'=S_3$, and $S_1'=\{s_3\}$.
We may assume that $s_1$, $s_2$, and $s_3$ are reflections along the planes $x=0$, $y=0$, and $z=0$, respectively, in $\RR^3$.
We use the standard basis $\{\ee_1,\ee_2,\ee_3\}$
to make the fundamental open chamber the first octant.
Then $H$ is the octahedral 3-graph on $\{X_0,X_1,Y_0,Y_1,Z_0,Z_1\}$ with edges $\{X_iY_jZ_k:i,j,k=0,1\}$,
where 
\begin{align*}
X_0&=\{(x,0,0):x>0\}, ~X_1=\{(x,0,0):x< 0\},\\
Y_0&=\{(0,y,0):y> 0\}, ~Y_1=\{(0,y,0):y< 0\},\\
Z_0&=\{(0,0,z):z> 0\}, ~Z_1=\{(0,0,z):z< 0\}.
\end{align*}
As noted in Example \ref{ex:Gowers_norms}, the $H$-norm is the Gowers' octahedral norm for $3$-graphs.
On the other hand, $H'$ is the 3-graph on eight vertices $\{X_{00},X_{01},X_{10},X_{11},Y_0,Y_1,Z_0,Z_1\}$ with eight edges $\{X_{ij}Y_iZ_k:i,j,k=0,1\}$.
Here $Y_i,Z_j$ are exactly the same as in $H$ since $S_2=S_2'$ and $S_3=S_3'$, 
and
\begin{align*}
&X_{00}=\{(x,y,0):x> 0,y> 0\},~X_{01}=\{(x,y,0):x<0,y> 0\},\\
&X_{10}=\{(x,y,0):x> 0,y< 0\},~X_{11}=\{(x,y,0):x< 0,y<0\}.
\end{align*}
Observe that the induced subgraphs of $H$ on $\{X_0,X_1,Y_i,Z_0,Z_1\}$, for $i=0,1$, and
the induced subgraphs of $H'$ on $\{X_{j0},X_{j1},Y_j,Z_0,Z_1\}$, for $j=0,1$, are 
all isomorphic to the kite-shaped 3-graph with $4$ edges, say $K$.
Moreover, $H$ is obtained by gluing two 
vertex-disjoint copies of $K$, identifying the copies of $X_0$, $X_1$, $Z_0$, and $Z_1$ in each graph.
We can build $H'$ similarly but now we only identify the copies of  
$Z_0$ and $Z_1$.
Therefore, setting
\begin{align*}
	g(x_0,x_1,z_0,z_1)=\int \prod_{i,j=0,1} f(x_i,y,z_j)dy,
\end{align*}
whose average is the $4^{th}$ power of the $K$-norm,
we have
\begin{align*}
	\norm{f}_{H}=\int g^2 dx_0dx_1dz_0dz_1
	\text{ and }\norm{f}_{H'}=\int \left(\int g dx_0dx_1\right)^2 dz_0dz_1.
\end{align*}
Thus, the Cauchy--Schwarz inequality implies that $\norm{f}_H\geq \norm{f}_{H'}$.
We shall generalise this idea later to prove Proposition \ref{prop:domination}.
\end{example}

\medskip

Another way to compare (hyper)graph norms is to ask that they be polynomially related, in the sense that if one norm is large, say at least $c$, then the other norm is at least $c^{k}$ for some appropriate $k$. 
Formally, we say that a (semi-)norm $\norm{\cdot}$ \emph{polynomially dominates} another (semi-)norm $\norm{\cdot}_0$ if, for any $|f|\leq 1$, 
\begin{align*}
	\norm{f}\leq c \text{ implies } 
	\norm{f}_{0}\leq c^\delta,
\end{align*}
where $\delta>0$ is a constant independent from $f$.
It is also possible for two norms to polynomially dominate each other, 
so in this case we say that the two norms are \emph{polynomially equivalent}.
For the rest of this subsection,
we shall state our results only for norming graphs, 
but analogous statements for weakly norming graphs hold if we replace $f$ by $|f|$ and graph norms by absolute graph norms.
Following Gowers' approach \cite{G01, G06, G07} to quasirandomness, 
but borrowing notation from graph limit theory \cite{L12,LSz06}, 
we define the cut-norm $\norm{f}_{\square}$ for a bounded measurable function $f$ on $[0,1]^2$ by
\begin{align*}
	\norm{f}_{\square}:=\sup_{u,v}\left\vert\int f(x,y)u(x)v(y)dxdy\right\vert,
\end{align*}
where the supremum is taken over all measurable $u,v:[0,1]\rightarrow[-1,1]$.
By applying inequality~\eqref{eqn:rainbow} with $C_4$, we see that
\begin{align}\label{eqn:uniform}
	\int f(x,y)u(x)v(y)dxdy \leq \norm{f}_{C_4}\norm{u}_{C_4}\norm{v}_{C_4}\leq \norm{f}_{C_4}
\end{align}
holds for all measurable $u,v:[0,1]\rightarrow[-1,1]$, and hence $\norm{f}_{\square}\leq\norm{f}_{C_4}$.
Conversely, if 
\begin{align*}
\norm{f}_{C_4}^4=\int f(x,y)f(x,y')f(x',y)f(x',y')> c^4,
\end{align*}
then there exist some $x^*,y^*\in[0,1]$
such that
\begin{align}\label{eqn:cutnorm_converse}
	\left\vert\int f(x,y)f(x,y^*)f(x^*,y)f(x^*,y^*)dxdy\right\vert \geq c^4.
\end{align}
Taking $u(x)=f(x^*,y^*)f(x,y^*)$ and $v(y)=f(x^*,y)$ gives $\norm{f}_{\square}\geq c^4$,
so the $C_4$-norm and the cut-norm are polynomially equivalent.
Similarly, using \eqref{eqn:rainbow} to deduce \eqref{eqn:uniform}
and averaging to obtain \eqref{eqn:cutnorm_converse},
it is easy to see that any $H$-norm is polynomially equivalent to the cut-norm
whenever $H$ is norming and contains a cycle.
Hence, if two norming graphs $H$ and $J$ contain cycles, then the $H$-norm and the $J$-norm are always polynomially equivalent.
Conversely, suppose $H$ is a norming graph isomorphic to a tree.
Then $H$ must be bi-regular by Theorem 2.10 in \cite{H10}
and, hence, is isomorphic to a star $K_{1,t}$.
By taking $f$ to be the balanced function of a non-quasirandom but regular graph,
one may obtain an example with small $K_{1,t}$-norm but large cut-norm.
Hence, star-norms are polynomially dominated by the cut-norm, but not vice versa.

Our aim now is to decide whether two different $H$-norms, where each $H$ is an $(S_1,\cdots,S_k;S,W)$-reflection hypergraph, are polynomially equivalent. The result below says that we may `forget' the angles between the reflection hyperplanes of a reflection group $W$
if $H$ is an $(S_1,\cdots,S_k;S,W)$-reflection hypergraph and still obtain a polynomially equivalent norm:
\begin{proposition}\label{prop:forget_angles}
 Let $W$ and $W'$ be two reflection groups with the same number of simple reflections, i.e.,
 $S=\{s_1,\cdots,s_n\}$ and $S'=\{s_1',\cdots,s_n'\}$ are the set of simple reflections of $W$ and $W'$, respectively.
 Suppose, for $S_1,\cdots,S_k\subset S$ and $S_1',\cdots,S_k'\subset S'$, that the
 indices covered by $S_i$ and $S_i'$ are always the same. That is, $\{j:s_j\in S_i\}=\{j:s_j'\in S_i'\}$ for each $i=1,2,\cdots,k$.
 Then the $(S_1,\cdots,S_k;S,W)$-reflection hypergraph $H$ and the $(S_1'\cdots,S_k';S',W')$-reflection hypergraph $H'$
 give polynomially equivalent hypergraph semi-norms,
 provided that both hypergraphs are norming.  
\end{proposition}
For example, the proposition above allows us to prove that 
if a reflection hypergraph $H$ consists of `tight triples', 
then the $H$-norm is polynomially equivalent to Gowers' octahedral norm for $3$-graphs.
\begin{example}\label{ex:S_4_hyper}
	In Example \ref{ex:S_4},
	the reflection group $W$ is isomorphic to the symmetric group on 4 elements
	and the set of simple reflections $S=\{t_{12},t_{23},t_{34}\}$,
	where $t_{ij}$ is a reflection swapping the $i$-th and $j$-th coordinates.
	Let $S_1=S\setminus\{t_{12}\}$, $S_2=S\setminus\{t_{23}\}$, $S_3=S\setminus\{t_{34}\}$,
	and $H$ be the $(S_1,S_2,S_3;S,W)$-reflection hypergraph.
	Observe that $H$ is a tripartite 3-graph isomorphic to the vertex-edge-face incidence $3$-graph of a tetrahedron,
	i.e., a vertex, an edge, and a face form an edge if they can be extended to a flag.
	The fact that $S_1\cap S_2\cap S_3=\emptyset$ implies that $H$ is norming and, 
	by Proposition \ref{prop:forget_angles} and Example~\ref{ex:Gowers_norms}, the $H$-norm is polynomially equivalent to Gowers' octahedral norm for $3$-graphs.	
\end{example}

The key to proving Proposition \ref{prop:forget_angles} is to show a polynomial equivalence between hypergraph norms and suitably generalised cut-norms.
For hypergraphs, Gowers~\cite{G06,G07} proved that there is a polynomial equivalence between the octahedral norms defined in Example~\ref{ex:Gowers_norms} and certain generalised cut-norms. Later, Conlon, H\`{a}n, Person, and Schacht \cite{CHPS12} showed that the norms discussed in Example~ \ref{ex:vertex-quasirandom} exist and are polynomially equivalent to certain weaker cut-norms. More recently, Reiher, R\"{o}dl, and Schacht \cite{RRS15} proposed further cut-norms lying between these two extremes. As we shall see, the result below gives a suitable $H$-norm which is polynomially equivalent to each of these cut-norms.

Let $\M=(M_1,M_2,\cdots,M_n)$ be an $n$-tuple of subsets of $[k]=\{1,2,\cdots,k\}$ and 
define the hypergraph cut-norm $\norm{\cdot}_{\square,\M}$ with respect to $\M$ by
	\begin{align*}
		\norm{f}_{\square,\M}:=\sup_{u_1,\cdots ,u_k}\left\vert\int f(x_1,\cdots,x_k)\prod_{i=1}^n u_i(x_{M_i})\right\vert,
	\end{align*}
where $f$ is a bounded measurable function on $[0,1]^k$, $x_I$ for $I\subset [k]$ is the vector $(x_i)_{i\in I}$ of variables with indices in $I$,
and the supremum is taken over all measurable functions $u_1,\cdots,u_k$ taking values between $[-1,1]$.
For example, $n=2$, $M_1=\{1\}$, and $M_2=\{2\}$ is the graph cut-norm.
The cut-norms described by Gowers \cite{G06,G07} 
correspond to the case when $n=k$ and $M_i=[k]\setminus\{i\}$.
In particular, for $n=3$, we have
\begin{align*}
	\norm{f}_{\square,\M}=\sup_{u,v,w:[0,1]^2\rightarrow[-1,1]}\left\vert\int  f(x,y,z)u(x,y)v(y,z)w(z,x)\right\vert.
\end{align*}
The lemma below states that whenever $H$ is the $(S_1,\cdots,S_k;S,W)$-reflection hypergraph, 
the $H$-norm is polynomially equivalent to the generalised cut-norm with respect to a suitably chosen $\M$.
\begin{lemma}
	Let $W$ be a reflection group with simple reflections $S = \{s_1, \cdots, s_n\}$
	and let $S_1,\cdots,S_k$ be subsets of $S$.
	Suppose $\M=(M_1,\cdots,M_n)$ is the $n$-tuple of subsets of $[k]$ 
	such that $M_i=\{j\in [k]:s_i\in S_j\}$ for $i=1,\cdots,n$.
	Then, whenever $H$ is an $(S_1,\cdots,S_k;S,W)$-reflection hypergraph that is norming, 
	the $H$-norm is polynomially equivalent to the cut-norm $\norm{\cdot}_{\square,\M}$.
\end{lemma}
\begin{proof}
Let $V_1V_2\cdots V_k$ be the edge of the $(S_1,S_2,\cdots,S_k;S,W)$-reflection hypergraph $H$ contained in the fundamental closed chamber.
Observe that $M_i$ is the collection of indices $j$ such that $V_j$ is fixed under the simple reflection $s_i$, since 
each collection of inactive mirrors $S_j$ defines the cone $V_j$. It follows that there is another edge of $H$, formed by reflecting 
the edge $V_1 V_2 \dots V_k$ in $s_i$, that contains the set $\{V_j: j \in M_i\}$. The remainder of the proof is now similar to the $C_4$ case discussed earlier.

To see that $\norm{f}_{\square,\M}\leq\norm{f}_H$, observe that, since $H$ is norming,
\begin{align*}
	\int f(x_1,\cdots,x_k)\prod_{i=1}^n u_i(x_{M_i})\leq\norm{f}_H\prod_{i=1}^n \norm{u_i}_H\leq \norm{f}_H
\end{align*}
whenever $u_i:[0,1]^{|M_i|}\rightarrow[-1,1]$.
To prove the converse, note that if $\norm{f}_H>c$, then
we may assign values to the variables corresponding to the vertices other than $V_1,V_2,\cdots,V_k$ to get
\begin{align*}
	\left\vert\int f(x_1,\cdots,x_k)\prod_{i=1}^n u_i(x_{M_i})\right\vert\geq c^{|E(H)|},
\end{align*}
as desired.
\end{proof}
Since the choice of $\M$ is independent of the group structure of $W$,
Proposition \ref{prop:forget_angles} follows as a corollary.
Note that each $M_i$ corresponds to the $i$-th column vector of the $k\times n$ 
incidence matrix between $k$ sets $S_1,S_2,\cdots,S_k$ and the simple reflections $s_1,s_2,\cdots,s_n$.
Thus, given a cut-norm $\norm{\cdot}_{\square,\M}$ with an ordered collection $\M$ of index subsets of $k$,
we may always construct a $k$-graph norm that is polynomially equivalent to it.
To give some examples, we will now describe hypergraph norms that are polynomially equivalent 
to the weak cut-norms introduced in \cite{RRS15}.
\begin{example}
In \cite{RRS15}, two different cut-norms for 3-graphs were studied.
For $3$-graphs, we always have $k=3$
and if we let $n=2$, $M_1=\{1\}$, and $M_2=\{2,3\}$, we get the cut-norm 
\begin{align*}
	\norm{f}_{\square,\M}=\sup_{u,v,w:[0,1]^2\rightarrow[-1,1]}\left\vert\int f(x,y,z)u(x)v(y,z)\right\vert.
\end{align*}
To ensure that $M_i=\{j:s_i\in S_j\}$ holds,
let $S_1=\{s_1\}$ and $S_2=S_3=\{s_2\}$.
Considering the simplest reflection group $W=\ang{s_1}\times\ang{s_2}$
generated by $s_1$ and $s_2$,
the $(S_1,S_2,S_3;S,W)$-reflection hypergraph is isomorphic to the $3$-graph $H$
with
\begin{align*}
	V(H)=\{x_0,x_1,y_0,y_1,z_0,z_1\}
	\text{ and }
	E(H)=\{x_iy_jz_j:i,j=0,1\}.
\end{align*}
The second cut-norm is the case when $n=2$, $M_1=\{1,2\}$, and $M_2=\{2,3\}$, i.e.,
\begin{align*}
	\norm{f}_{\square,\M}=\sup_{u,v,w:[0,1]^2\rightarrow[-1,1]}\left\vert\int f(x,y,z)u(x,y)v(y,z)\right\vert.
\end{align*}
Then we have $S_1=\{s_1\},S_2=\{s_1,s_2\}$, and $S_3=\{s_2\}$.
With the same group $W=\ang{s_1}\times\ang{s_2}$,
the $(S_1,S_2,S_3;S,W)$-reflection hypergraph is isomorphic to the kite-shaped 3-graph with $5$ vertices
denoted by $K$ in Example \ref{ex:domination}.
\end{example}

By combining Proposition \ref{prop:forget_angles} and Proposition \ref{prop:domination}, 
we also have the following weak domination result between hypergraph norms with the same cut involution group.
\begin{corollary}
	Let $W$ be a reflection group with $S=\{s_1,\cdots,s_n\}$ the set of simple reflections
	and let $S_1,\cdots,S_k,S_1',\cdots,S_k'$ be subsets of $S$ such that 
	$S_i'\subset S_i$ for each $i=1,\cdots,k$.
	Suppose $H$ and $H'$ are the $(S_1,\cdots,S_k;S,W)$-reflection hypergraph and the $(S_1',\cdots,S_k';S,W)$-reflection hypergraph, respectively, and suppose that they are both norming.
	Then the $H$-norm polynomially dominates the $H'$-norm.
\end{corollary}

\medskip

Getting back to the proof of Proposition \ref{prop:domination},
we first describe the structure of the given hypergraphs $H$ and $H'$.
Let $W\subset\mathbf{GL}(n,\RR)$ be a reflection group with simple reflections $S=\{s_1,s_2,\cdots,s_n\}$
such that each $s_i$ is the reflection along the hyperplane
\begin{align*}
 H_i:=\{(x_1,x_2,\cdots,x_n)\in\RR^n: x_i=0\},
\end{align*}
i.e., $s_i$ flips the sign of the $i$-th coordinate.
Then the $(S_1,\cdots,S_k;S,W)$-reflection hypergraph $H$ is the $k$-partite graph with $k$-partition $\A_1\cup\cdots\cup\A_k$,
where $\A_i$ consists of the $2^{n-|S_i|}$ vertices of the form 
\begin{align*}
	v_{i}(\xi):=\{(x_1,\cdots,x_n):\xi_j x_j> 0
	~\mbox{for all } j\notin S_i,~ x_j=0 ~\mbox{for all } j\in S_i\},
\end{align*}
where $~\xi=(\xi_1,\cdots,\xi_n)\in\{\pm 1\}^n$.
Here each $\xi\in\{\pm 1\}^n$ represents the closed chamber
\begin{align*}
\{(x_1,\cdots,x_n):\xi_j x_j\geq 0
	~\mbox{for all } j=1,2,\cdots,n\},
\end{align*}
so a $k$-tuple of vertices in $\A_1\times\A_2\times\cdots\times\A_k$
is an edge of $H$ 
if and only if it can be written as $(v_1(\xi),v_2(\xi),\cdots,v_k(\xi))$ with the same $\xi\in\{\pm 1\}^n$.
Let $H^+$ be the subgraph of $H$ induced on cones in the half-space $\{(x_1,\cdots,x_n):x_n\geq 0\}$ induced by the hyperplane $H_n$.
Using the fact that $s_n$ is a cut involution of $H$, we have
\begin{align}\label{eqn:bond}
	\norm{f}_H^{|E(H)|}=\int\left(\int g ~dy_{V(H^+)\setminus F(H)}\right)^2 dy_{F(H)},
\end{align} 
where $g=\prod_{(u_1,\cdots,u_k)\in E(H^+)}f(y_{u_1},\cdots,y_{u_k})$, $y_I=(y_i)_{i\in I}$ for $I\subset V(H)$, 
and $F(H)$ is the set of vertices $\{v_j(\xi)\in V(H): s_n\in S_j, \xi\in\{\pm 1\}^n\}$
that are fixed by $s_n$.
Furthermore, $H^+$ is isomorphic to the $(S_1^+,\cdots,S_k^+;S,W)$-reflection hypergraph with $S_i^+=S_i\cup\{s_n\}$ for all $i$, 
since adding $s_n$ to all $S_i$ makes the sign of $x_n$ fixed.

\begin{proof}[Proof of Proposition \ref{prop:domination}]
	To show that $H$ is always norming, we may assume that
	there exists a simple reflection $s\in\bigcap_{i=1}^k S_i$.
	Observe that the $(S_1^-,\cdots,S_k^-;S,W)$-reflection hypergraph $H^-$, where $S_i^-=S_i\setminus \{s\}$ for all $i$, 
	consists of two vertex-disjoint copies of $H$ mapped to each other by $s$.
	Then the $H^-$-norm takes exactly the same value as the $H$-norm and, hence, by induction on $|\bigcap_{i=1}^k S_i|$, we are done.
	
 	To prove that the $H$-norm dominates the $H'$-norm,
 	we may assume $S_i'=S_i\setminus\{s_n\}$, $s_n\in S_i$, and $S_j'=S_j$ for $j\neq i$.
	Both $H$ and $H'$ have the same induced subgraph on cones in the half-space $\{(x_1,\cdots,x_n):x_n\geq 0\}$,
	because it must be isomorphic to the $(S_1^+,\cdots,S_k^+;S,W)$-reflection hypergraph with $S_j^+=S_j\cup\{s_n\}$ for all $j$.
	When writing $\norm{f}_{H'}$ as in \eqref{eqn:bond},
	the fixed point set $F(H')$ under $s_n$ becomes
	\begin{align*}
	F(H')=\{v_j(\xi)\in V(H'): s_n\in S_j', \xi\in\{\pm 1\}^n\}=F(H)\setminus \A_i,
	\end{align*}
	since $s_n$ is removed from $S_i$ .
	Therefore, the $H'$-norm can be expressed, with the same $g$ as in \eqref{eqn:bond}, by
	\begin{align*}
		\norm{f}_{H'}^{|E(H')|}=\int\left(\int g ~dy_{V(H^+)\setminus F(H)}dy_{\A_i}\right)^2 dy_{F(H)\setminus \A_i},
	\end{align*}
	and the Cauchy--Schwarz inequality implies that this is less than or equal to \eqref{eqn:bond}.
\end{proof}

\subsection{Applications to Sidorenko's conjecture} \label{sec:Sid}

We will now discuss how to apply our results to Sidorenko's conjecture. As mentioned in the introduction, there is the obvious fact that any weakly norming graph satisfies Sidorenko's conjecture. However, one can go beyond this by applying the entropy techniques that have been developed in recent years~\cite{KR11, LSz12, KLL14, Sz15, CKLL15}. The key facts we will use about entropy are contained in the following lemma, though we refer the reader to \cite{CKLL15} for further information on entropy, conditional entropy, and their use in relation to Sidorenko's conjecture. Throughout this subsection, logarithms will be understood to be base 2.

\begin{lemma}\label{lem:entropy}
Let $X$, $Y$, and $Z$ be random variables and suppose that $X$ takes values in 
a set $S$, $\HH(X)$ is the entropy of $X$, and $\HH(X|Y)$ is the conditional entropy of $X$ given $Y$. Then 
\begin{enumerate}[(i)]
	\item $\HH(X)\leq\log|S|$,
	\item $\HH(X|Y,Z)=\HH(X|Z)$ if $X$ and $Y$ are conditionally independent given $Z$.
\end{enumerate}
\end{lemma}
As in \cite{KLL14,CKLL15}, we say that a bipartite graph $H$ has Sidorenko's property if $H$ satisfies \eqref{eqn:Sidorenko}
for all non-negative symmetric functions $f$, i.e., Sidorenko's conjecture holds for $H$.

\medskip

To motivate what follows, suppose we wish to show that the graph obtained by gluing two copies $H_1$ and $H_2$ of the 1-subdivision of $K_4$ along induced 6-cycles has Sidorenko's property.
Formally, there is a vertex set $J= V(H_1)\cap V(H_2)$ on which the induced subgraphs of both $H_1$ and $H_2$ are isomorphic to a 6-cycle 
and we are interested in proving that the graph $H$ on vertex set $V(H_1)\cup V(H_2)$ with edges $E(H_1)\cup E(H_2)$ has Sidorenko's property.
As in \cite{CKLL15} and \cite{Sz15},
we generate a random homomorphism in $\Hom(H,G)$ and analyse its entropy to get a lower bound for $|\Hom(H,G)|$.
To this end, we consider the following way of generating a random element in $\Hom(H,G)$:
\begin{enumerate}[(i)]
	\item Take a uniform random homomorphic copy of $H_1$, i.e., choose an element in $\Hom(H_1,G)$ uniformly at random.
	\item Choose another homomorphic copy of $H_2$ that extends the 6-cycle on $J$ uniformly at random.
\end{enumerate}
This algorithm always gives a random homomorphism $\ww$ in $\Hom(H,G)$, so it remains to analyse the entropy.
We regard $\ww$ as a random vector indexed by $V(H)$ and, for $I\subset V(H)$, write $\ww_I$ for the random vector $(\ww(v))_{v\in I}\in V(G)^I$.
Let $U=V(H_1)\setminus V(J)$ and $W=V(H_2)\setminus V(J)$.
Then
\begin{align}\label{eqn:entropy_motivate}
\nonumber
	\log|\Hom(H,G)|&\geq \HH(\ww)=\HH(\ww_{U},\ww_J,\ww_W)\\ \nonumber
					&=\HH(\ww_W|\ww_U,\ww_J)+\HH(\ww_U,\ww_J)\\ \nonumber
					&=\HH(\ww_W|\ww_J)+\log|\Hom(H_1,G)|\\ \nonumber
					&=\HH(\ww_W,\ww_J)-\HH(\ww_J)+\log|\Hom(H_1,G)|\\
					&=2\log|\Hom(H_1,G)|-\HH(\ww_J).
\end{align}
If $J$ were an independent set, we could upper bound $|\Hom(J, G)|$ by $|G|^{|J|}$ and use Lemma \ref{lem:entropy}~(i) to derive an upper bound on $\HH(\ww_J)$. This is one of the main tricks used in \cite{CKLL15} and \cite{Sz15}, an elaboration of which also allows one to handle the case where the graph induced on $J$ is a tree. 
In our case, the graph induced on $J$ is a $6$-cycle,
so we have $\HH(\ww_J)\leq\log|\Hom(C_6,G)|$.
Since $H_1$ is a weakly norming graph, the monotonicity property~\eqref{eqn:mono}
implies that the homomorphism densities of $H_1$ and $C_6$ in $G$ satisfy $t_{H_1}(G)^{1/|E(H_1)|} \geq t_{C_6}(G)^{1/|E(C_6|}$. This in turn implies that
\begin{align*}
	|\Hom(H_1,G)|^{1/2}|V(G)|\geq |\Hom(C_6,G)|.
\end{align*}
Plugging this into \eqref{eqn:entropy_motivate}, we have
\begin{align*}
	\log|\Hom(H,G)|&\geq \frac{3}{2}\log|\Hom(H_1,G)|-\log|V(G)|\\
					&\geq 18\log\left(\frac{|\Hom(K_2,G)|}{|V(G)|^2}\right)+14\log|V(G)|,
\end{align*}
where the last inequality follows from the fact that $H_1$ has Sidorenko's property.
Therefore, $H$ has Sidorenko's property.

More generally, we may iterate this gluing process in a tree-like way.
To describe the resulting graph, we use the notion of tree decomposition introduced by Halin~\cite{Hal76}
and developed by Robertson and Seymour \cite{RS84}.
A \emph{tree decomposition} of a graph $H$ is a pair $(\mathcal{F}, \TT)$ consisting of a family $\mathcal{F}$  of vertex subsets of $H$ and  a tree $\TT$ on 
vertex set $\mathcal{F}$ satisfying
\begin{enumerate}
\item $\bigcup_{X\in\mathcal{F}}X=V(H)$, 
\item for each $\{v,w\} \in E(H)$, there exists a set $X \in \mathcal{F}$ such that
$v, w \in X$, and
\item for $X,Y,Z\in \mathcal{F}$, $X\cap Y\subseteq Z$ 
whenever $Z$ lies on the path from $X$ to $Y$ in $\TT$.
\end{enumerate}
Now, given a weakly norming graph $N$, an \emph{$N$-decomposition} of a graph $H$ is a tree decomposition $(\mathcal{F},\TT)$ of $H$ 
satisfying the following two extra conditions:
\begin{enumerate}
\item 
The induced subgraphs $H[X]$, $X \in \FF$, are each isomorphic to $N$.
\item 
For every pair $X,Y\in \FF$ which are adjacent in $\TT$, there is an isomorphism between the two copies of $N$ that fixes $X \cap Y$.
\end{enumerate}
We say that a bipartite graph is \emph{$N$-decomposable} if it allows a $N$-decomposition,
i.e., it can be obtained by gluing copies of the weakly norming graph $N$ in a tree-like way.
The main theorem in this subsection is as follows.
\begin{theorem}\label{thm:tree_gluing}
	If $N$ is weakly norming, any $N$-decomposable graph $H$ has Sidorenko's property.
\end{theorem}
\begin{proof}
	Let $(\FF,\TT)$ be an $N$-decomposition of $H$ and let $G$ be the target graph in which we wish to embed $H$.
	The randomised algorithm for generating a copy of $H$ is a straightforward generalisation of that discussed in the example above:
	pick a root $R\in \FF$ and choose a uniform random homomorphism from $\Hom(H[R],G)$
	and, for each child $X$ of $R$, choose a uniform random (homomorphic) copy of $H[X]$ extending 
	the embedded copy of $H[R\cap X]$.
	Repeating this process, we obtain a random homomorphic copy $\ww$ of $H$.
	Following the proof of Theorem 1.2 in  \cite{CKLL15}, we arrive at the identity 
	 \begin{align*}
	 	\HH(\ww)=\sum_{X\in\FF}\HH(\ww_X)-\sum_{XY\in E(\TT)}\HH(\ww_{X\cap Y}).
	 \end{align*}
	 Let $K_{XY}$ be the induced subgraph $H[X\cap Y]$ for $XY\in E(\TT)$.
	 We may bound $\HH(\ww)$ and $\HH(\ww_{X\cap Y})$ from above by $\log|\Hom(H,G)|$
	 and $\log|\Hom(K_{XY},G)|$, respectively. Moreover, since $\ww_X$ is distributed uniformly in
	 $\Hom(N, G)$, this gives 
	 \begin{align}\label{eqn:norming_tree}
	 	\log t_H(G)\geq |\FF|\log t_N(G)-\sum_{XY\in E(\TT)} \log t_{K_{XY}}(G),
	 \end{align}
	 after rescaling the equation by subtracting $|V(H)|\log |V(G)|$ on each side.
	 Since $N$ is weakly norming, 
	 the monotonicity property \eqref{eqn:mono} implies that
	 \begin{align*}
	 	t_{K_{XY}}(G)\leq t_N(G)^{|E(K_{XY})|/|E(N)|}.
	\end{align*}
	Plugging this bound into \eqref{eqn:norming_tree}, it follows that
	\begin{align*}
		\log t_H(G)\geq \frac{1}{|E(N)|}\left(|\FF||E(N)|-\sum_{XY\in E(\TT)}|E(K_{XY})|\right)\log t_N(G).
	\end{align*}
	Note that $|\FF||E(N)|-\sum_{XY\in E(\TT)}|E(K_{XY})|$ is exactly the number of edges in $H$
	and, hence, we have the bound $t_H(G)\geq t_N(G)^{|E(H)|/|E(N)|}$.
	As $N$ has Sidorenko's property, we conclude that $H$ does also.
\end{proof}

The theorem above is intended as an example of how norming graphs may be used to build up new examples of graphs satisfying Sidorenko's property. We expect that additional refinements of the entropy method can be used to broaden this class further, but we have not attempted a comprehensive treatment here. 

\section{Concluding remarks}
\textbf{Characterising weakly norming graphs.}
Although Theorem \ref{thm:parabolic} gives a fairly large class of weakly norming graphs,
it is still an open problem to characterise them all. For our arguments to work on a particular graph $H$,
it was necessary that $H$ be edge-transitive under its cut involution group. We also insisted on the
existence of a percolating sequence in $H$, but we suspect that this condition may not be needed.

\begin{conjecture}\label{conj:edge-trans}
	A bipartite graph $H$ is weakly norming if it is edge-transitive under its cut involution group.
\end{conjecture}

Underlying this conjecture is a further suspicion that reflection graphs may constitute the entire class of graphs which 
are edge-transitive under their cut involution groups. Were this indeed the case, 
Conjecture \ref{conj:edge-trans} would follow immediately from Theorem \ref{thm:parabolic}.  

In \cite{H10}, Hatami showed that weakly norming graphs must be balanced and bi-regular
and that the class of weakly norming graphs is closed under taking tensor products. 
There are analogues of these results for reflection graphs: 
any reflection graph must be bi-regular
since it is vertex-transitive on each side, 
while the tensor product of two reflection graphs is again a reflection graph.
It is not too difficult to prove that the same properties hold for graphs that are edge-transitive under the cut involution group.
These coincidences suggest that the converse to Conjecture \ref{conj:edge-trans} may also be true, though this is significantly more tentative than Conjecture~\ref{conj:edge-trans}.

\begin{conjecture}\label{conj:converse}
	A bipartite graph $H$ is weakly norming only if it is edge-transitive under its cut involution group.
\end{conjecture}

It would already be very interesting to show that weakly norming graphs are necessarily edge-transitive.

\medskip

\noindent
\textbf{Norming graphs versus weakly norming graphs.}
When can we guarantee that a weakly norming graph is also norming?
Theorem \ref{thm:norming} gives a partial answer to this natural question.
However, we again suspect, in analogy with Conjectures \ref{conj:edge-trans} and \ref{conj:converse}, that edge-transitivity under the stable involution group is a necessary and sufficient property.

\begin{conjecture}\label{conj:converse_norm}
	A bipartite graph $H$ is norming if and only if it is edge-transitive under its stable involution group.
\end{conjecture}

This is reminiscent of another conjecture, proposed in \cite{CCHLL12}, about determining those graphs $H$ for which $t_H(f)$ is always non-negative. It might be interesting to investigate the connection.

\begin{conjecture}[Positive graph conjecture \cite{CCHLL12}]
	$t_H(f)$ is non-negative for all $f:[0,1]^2\rightarrow\RR$ if and only if there is a stable involution of $H$.
\end{conjecture}

\medskip

\noindent
\textbf{Complex-valued functions.}
Our results on norming (hyper)graphs, Theorem \ref{thm:norming} and Theorem \ref{thm:hypergraph}, can also be used to define 
(hyper)graph norms for complex-valued functions.
Suppose $H$ is the $(S_1,\cdots,S_k;S,W)$-reflection hypergraph, 
where $W$ is a reflection group and $\bigcap_{i=1}^k S_i=\emptyset$.
Note that each reflection in $W$ acts as a stable involution
and there exists a one-to-one correspondence between the edges of $H$ and the set $\C$ of open chambers, 
since each edge is contained in exactly one closed chamber.
By the one-to-one correspondence between $W$ and $\C$, there also exists a one-to-one correspondence $\xi:W\rightarrow E(H)$ between $W$ and $E(H)$.
Let $W_+:=\{w\in W:\ell(w)\text{ is even}\}$ and $W_-:=\{w\in W:\ell(w)\text{ is odd}\}$.
We may write $W_+$ and $W_-$ as the inverse images of $1$ and $-1$ under the group homomorphism $w\mapsto \det(w)=(-1)^{\ell(w)}$ from $W$ to the multiplicative group $\{\pm 1\}$.
It follows that $W_+$ and $W_-$ partition $W$ into parts of equal size
and each reflection $t\in W$ maps $w\in W_+$ to $tw\in W_-$ and vice versa.
For a complex-valued function $f$ on $[0,1]^k$, define
\begin{align*}
	\norm{f}_H:=
	\left\vert\int \prod_{e\in\xi(W_+)}f(x_e)
	\prod_{e'\in\xi(W_-)}\overbar{f(x_{e'})}
	~\right\vert^{1/|E(H)|},
\end{align*}
where $x_e=(x_{u_1},x_{u_2},\cdots,x_{u_k})$ for $e=(u_1,u_2,\cdots,u_k)$.
This is well-defined, even without taking the absolute value of the integral above,
since 
\begin{align*}
	\int \prod_{e\in\xi(W_+)}f(x_e)
	\prod_{e'\in\xi(W_-)}\overbar{f(x_{e'})}
	=\int g\overbar{g},
\end{align*}
where 
\begin{align*}
	g(x_{F_\phi})=\int \prod_{e\in\xi(W_+)\cap (L_\phi\cup F_\phi)^k}f(x_e)
\prod_{e'\in\xi(W_-)\cap (L_\phi\cup F_\phi)^k}\overbar{f(x_{e'})}~dx_{L_\phi}
\end{align*}
for a stable involution $\phi$ that corresponds to a reflection in $W$.
When $W=\ang{s_1}\times\cdots\times\ang{s_n}$, $s_i\in S$, this gives complex-valued versions of Gowers' octahedral norms, but it also
allows some more exotic examples. For example, the hypergraph described in Example \ref{ex:S_4_hyper} can now easily be modified to define a norm
on the vector space of complex-valued functions $f: [0,1]^3 \rightarrow \mathbb{C}$.

\medskip

\vspace*{5mm}

\textbf{Acknowledgements.}
Part of this work was carried out while the authors participated in the LMS\red{-}CMI Research School on Regularity and Analytic Methods in Combinatorics at the University of Warwick and also while the second author was visiting KIAS.
The second author would like to thank Seung Jin Lee for suggestions of references and helpful discussions on algebraic combinatorics. We would also like to thank Alexander Sidorenko for some helpful remarks on an earlier version of this paper.

\bibliographystyle{abbrv}
\bibliography{references}

\end{document}